\documentclass[leqno,11pt]{article}

\usepackage{amssymb,hyperref,amsthm,hyperref}
\usepackage{euscript}
\usepackage[dvips]{graphicx}
\usepackage{flafter}
\usepackage{pstricks}
\usepackage[latin1]{inputenc}
\usepackage{amsmath}
\usepackage{amsfonts}
\usepackage{pstricks-add}
\usepackage{dsfont}
\usepackage{bm}
\usepackage{tikz}

\hypersetup{
    linktoc=page,
    linkcolor=red,          % color of internal links
    citecolor=blue,        % color of links to bibliography
    filecolor=blue,      % color of file links
    urlcolor=cyan,
    colorlinks=true           % color of external links
}

\usepackage[refpage]{nomencl}
\makenomenclature

\usepackage{makeidx}
\makeindex

\usepackage{mathrsfs} %%% Nice caligraphic fonts
\evensidemargin\oddsidemargin
\begin{document}

%%
%%%%%%%%%%%%%%%%%%%%%%%%%%% Equation numberings
\newcommand{\eqnsection}{
\renewcommand{\theequation}{\thesection.\arabic{equation}}
   \makeatletter
   \csname  @addtoreset\endcsname{equation}{section}
   \makeatother}
\eqnsection
%%%%%%%%%%%%%%%%%%%%%%%%%%%

%%%%%%%%%%%%%% Bbb characters
%%%%%%%%%%%%%% Real numbers
\def\r{{\mathbb R}}
\def\u{ \mathtt u}
%%%%%%%%%%%%%% Expectation
\def\e{{\mathbb E}}
%%%%%%%%%%%%%% Probability
\def\p{{\mathbb P}}
%%%%%%%%%%%%%% Law of environment
\def\P{{\mathbb P}}
\def\E{{\mathbb E}}
\def\Q{{\bf Q}}
%%%%%%%%%%%%%% Integers
\def\rc{]\!]}
\def\lc{[\![}
\def\z{{\mathbb Z}}
%%%%%%%%%%%%%% Natural numbers
\def\N{{\mathbb N}}
%%%%%%%%%%%%%% Tree
\def\T{{\mathbb T}}
%%%%%%%%%%%%%% Galton-Watson tree
\def\G{{\mathbb G}}

\def\C{{\mathbb C}}
%%%%%%%%%%%%%% "Low"
\def\L{{\mathbb L}}
%%%%%%%%%%%%%% "indicatrice"
\def\1{{\mathds{1}}}
%%%%%%%%%%%%%%
\def\deg{\chi}

\def\d{\mathtt{d}}
%%%%%%%%%%%%%%
\def\ttheta{{\bm \theta}}
%%%%%%%%%%%%%%
%%%%%%%%%%%%%% "vecteur t"
\def\t{{\bf{t}}}
\def\a{{\bf{a}}}
%%%%%%%%%%%%%%
\def\deg{\chi}
\def\B{\mathbb{B}}
\def\V{{\mathtt{V}}} 
%%%%%%%%%%%%%%
%%%%%%%%%%%%%%%% Special symbols
%%%%%%%%%%%%%% Exponential
\def\ee{e}
%%%%%%%%%%%%%% Differentiation
\def\d{\, \mathrm{d}}
%%%%%%%%%%%%%% Survival
\def\S{\mathscr{S}}
%%%%%%%%%%%%%% Binary search
\def\bs{{\tt bs}}
%%%%%%%%%%%%%%
\def\bbeta{{\bm \beta}}
%%%%%%%%%%%%%%
\def\M{\mathtt{M}}

\def\ttheta{{\bm \theta}}
%%%%%%%%%%%%%%
%%%%%%%%%%%%%%%%%%%%%%%%%%%%%%%%%%%%%%%%%%%%
%%%%%%%%%%%%%%%%%%%%%%%%%%%%%%%%%%%%%%%%%%%%
%%%%%%%%%%%%%%%%%%%%%%%%%%%%%%%%%%%%%%%%%%%%

\newtheorem{theorem}{Theorem}[section]

\newtheorem{definition}[theorem]{Definition}
\newtheorem{Lemma}[theorem]{Lemma}
\newtheorem{Proposition}[theorem]{Proposition}
\newtheorem{Remark}[theorem]{Remark}
\newtheorem{corollary}[theorem]{Corollary}
\newtheorem{conjecture}[theorem]{Conjecture}
%%%%%%%%%%%%%%%%%%%%%%%%%%%%%%%%%%%%%%%%%%%%
%%%%%%%%%%%%%%%%%%%%%%%%%%%%%%%%%%%%%%%%%%%%
%%%%%%%%%%%%%%%%%%%%%%%%%%%%%%%%%%%%%%%%%%%%

%%%%%%%%%%%%%% Beginning of the text

\title{The tail distribution of the Derivative martingale and the global minimum of the branching random walk}
\date{\vspace{-5ex}}
 
 \author{ Thomas Madaule \footnote{Universit\'e Paul Sabatier, IMT, Toulouse, France} } 
\date{\today}
 \maketitle

\begin{abstract}
In a seminal paper Biggins and Kyprianou \cite{BKy04} proved the existence of a non degenerate limit for the {\it Derivative martingale} of the branching random walk. As shown in \cite{Aid11} and \cite{Mad11}, this is an object of central importance in the study of the extremes of the branching random walk. In this note we investigate the tail distribution of the limit of the {\it Derivative Martingale} by mean of the study of the global minimum of the branching random walk. This new approach leads us to extend the results of \cite{Gui90} and \cite{Bur09} under slighter assumptions.  
\end{abstract}
 
\noindent{\bf Key words or phrases:}  branching random walk.

\noindent{\bf MSC 2000 subject classifications:  60J80, 60G57, 60G50,  60G17.}

{\leftskip=2truecm \rightskip=2truecm \baselineskip=15pt \small

} %%%%%% End of narrower

%%%%%%%%%%%%%%%%%%%%%%%%%%%%%%%%%%%%%%%%%%%%%%%%%%%%%%%%%%%%%%%%%%%%%%%%%%%%%%%%
\section{Introduction}
%%%%%%%%%%%%%%%%%%%%%%%%%%%%%%%%%%%%%%%%%%%%%%%%%%%%%%%%%%%%%%%%%%%%%%%%%%%%%%%%

 We consider a real-valued branching random walk{\string:} Initially,  a single particle denoted $\emptyset$ sits at the origin.  Its children together with their  displacements,  form a point process $\Theta$ on $\r$ and  the first generation of the branching random walk.  These children have children of their own which form the second generation, and behave --relatively to their respective positions at birth-- like independent copies of the same point process $\Theta$.  And so on.
 
 Let $\mathbb{T}$ be the genealogical tree of the particles in the branching random walk. Plainly, $\mathbb{T}$ is a Galton-Watson tree.  We write $|z|=n$ if a particle  $z$  is in the $n$-th generation, and denote its position by $V(z)$ ($V(\emptyset)=0$). The collection of positions $(V(z),z\in \mathbb{T})$ is our branching random walk.

We assume throughout the paper the following conditions: the distribution of $\Theta$ is non-lattice and
 \begin{eqnarray}
 \label{chap3criticalcondition1}
 \E\Big(\underset{|x|=1}{\sum}  \ee^{-V(x)} \Big)=1,&   &\E\Big(   \underset{|x|=1}{\sum}1  \Big) >1 ,\quad \text{  and }
 \\
 \label{chap3criticalcondition2} \E\Big(\underset{|x|=1}{\sum}V(x)\ee^{-V(x)} \Big)=0,&   \qquad &\sigma^2:=\E \Big(\underset{|x|=1}{\sum}V(x)^2\ee^{-V(x)} \Big)<\infty.
 \end{eqnarray}
The branching random walk is then said to be in the boundary case (Biggins and Kyprianou \cite{Bky05}). We refer to (the ArXiv version of) \cite{Jaf09} for detailed discussions on the nature of the assumption (\ref{chap3criticalcondition1}) and (\ref{chap3criticalcondition2})).

% 
% Let $\Phi(\beta):=\log \E \Big(\underset{|x|=1}{\sum}\ee^{-\beta V(x)} \Big)  \in (-\infty, +\infty],\, \beta \in \r $ and let
% \begin{equation}
% W_{\beta,n}:=\underset{|x|=n}{\sum}\ee^{-\beta V(x)-\Phi(\beta)n},\quad \beta\in \r,
% \end{equation} 
% which can be viewed as the normalized partition function of a directed polymer on trees, see the forthcoming (\ref{chap3defmubeta}). 
% In the literature $W_{1,n}= {\sum}_{|x|=n}\ee^{- V(x)} $ (note that $\Phi(1)=0$) is called the {\it critical additive martingale} associated with the branching random walk. For notational simplification, we write $W_{1,n}=W_n$ for any $n \geq 0$ [$W_0:=1$].   

Let us define respectively the critical  {\it additive and derivative Martingale}:
 \begin{align}
 W_n:= \sum_{|z|=n}\ee^{-V(z)},\qquad  D_n:= \sum_{|z|=n}V(z)\ee^{-V(z)},\qquad n\geq 0
 \end{align}
Defining $X:=\underset{|x|=1}{\sum}\ee^{-V(x)}$ 	and $ \tilde{X}:=\underset{|x|=1}{\sum}\max\{0,V(x)\}\ee^{-V(x)}$, it is well known, see Lyons \cite{Lyo97} and  Biggins and Kyprianou \cite{BKy04} that under  
 \begin{eqnarray}
 \label{chap31.5}
 &&\E(X(\max(0,\log X)^2)<\infty, \qquad  \E(\tilde{X}\max (0,\log \tilde{X}))<\infty,
 \end{eqnarray}
we have,
 \begin{equation}
\label{immedia}
\lim_{n\to \infty} W_n=0,\quad\qquad   \underset{n\to\infty}{\lim}D_n = D_\infty,\quad \P- \text{a.s}..
\end{equation}
Moreover the random variable $D_\infty$ is strictly positive on the set of non-extinction. Recently Chen \cite{Che15} proved, for a branching random walk in the boundary case, that the convergence of $(D_n)_{n\in \N}$ implies assumption (\ref{chap31.5}).  Both martingales are fundamental objects which have attracted many works this last decade. For instance they play a crucial role in the study of the extremes of the branching random walk, we cite a remarkable result due to Aïdékon:
\begin{theorem}[Aïdékon, \cite{Aid11}]
	\label{AIDi}
(\ref{chap3criticalcondition1}), (\ref{chap3criticalcondition2}) and  (\ref{chap31.5}). There exists a constant $C^*>0$ such that
\begin{equation}
  \lim_{n\to\infty} \P\left( \min_{|z|=n}V(z)- \frac{3}{2}\log n \geq x \right)=\E\left( \ee^{-C^* D_\infty \ee^x} \right).
\end{equation}
\end{theorem}

\paragraph{}In this paper we study, under slightly stronger assumption than (\ref{chap31.5}) the tail distribution of $D_\infty$. To our knowledge this question has only be tackled by mean of the study of the {\it smoothing transform}, i.e the study of the random variables $Z$ solution of
\begin{equation}
\label{Smooth}
Z \overset{(d)}{=} \sum_{|z|=1}\ee^{-\beta V(z)-\Phi(\beta)}Z^{(z)}, \qquad \beta>0,
\end{equation}
with $\Phi(\beta):= \log \E\big( \sum_{|z|=1} \ee^{-V(z)}   \big)$ and $(Z^{(z)})_{|z|=1}$ are independent copies of $Z$. It is easy to check that $D_\infty$ satisfies (\ref{Smooth}) with $\beta=1$ ($W_\infty=0$ is a degenerated solution of (\ref{Smooth})). This approach, first initiated by \cite{DLi83}, was very successful. Concerning the tail distribution of $D_\infty$ the most general result is due to Buraczewski \cite{Bur09}, he proved the following theorem:
\begin{theorem}[Buraczewski]
	\label{BBURA}
Consider a branching random walk satisfying (\ref{chap3criticalcondition1}), (\ref{chap3criticalcondition2}) and
\begin{equation}
\label{condBur}
\E\Big( \sum_{|z|=1}\ee^{-(1-\delta_1)V(z)}\Big) +\E\Big((\sum_{|z|=1}\ee^{-V(z)})^{1+\delta_2}\Big)<+\infty.
\end{equation}
with some constant $\delta_1,\delta_2>0$. Then there exists a strictly positive constant $C_0$ such that any non-negative solution $Z$ of (\ref{Smooth}), satisfies
\begin{equation}
\label{bubur}
\lim_{x\to\infty}x\P(Z>x)=C_0.
\end{equation}
\end{theorem} 
Theorem \ref{BBURA} extends a result of Guivarch \cite{Gui90} via similar techniques.
\begin{Remark}
	A key result to prove Theorem \ref{AIDi} consists to establish the tail distribution of $M_n:= \min_{|z|=n}V(x)$, i.e $\lim_{x\to \infty}\lim_{n\to\infty} \frac{\ee^x}{x}\P(M_n-\frac{3}{2}\ln n \leq -x)=C^*$. With standard manipulations one can deduce that
	%	{\yellow \begin{eqnarray}
	%		\E(\ee^{-q D})= 1 -q (\ln C^* + \ln \frac{1}{q}), \quad q\to 0
	%		\end{eqnarray}
	%		and then}
	\begin{equation}
		\int_0^x \P\left( D_\infty \geq u\right) du \sim \ln x,\qquad \text{when } x \to \infty.
	\end{equation}
	However it does not suffice to obtain the tail distribution of $D_\infty$. Indeed we are just in the situation where the monotone density theorem does not apply. 
\end{Remark}

The purpose of the present note is to prove the following  theorem  
\begin{theorem}
	\label{ourasi}
	Assume (\ref{chap3criticalcondition1}), (\ref{chap3criticalcondition2}) and 
	\begin{align}
	\label{ourAssum}
	\E((\tilde{X}+X)\max (0,\ln \tilde{X}+X)^5)<\infty.
	\end{align} 
There exists $c_{D_\infty}>0 $ such that
	\begin{equation}
	\label{TiTail}
	\lim_{x\to \infty} x \P(D_\infty \geq x) =c_{D_\infty}.
	\end{equation}
Moreover $c_{D_\infty}:= c_\M \times  \E\left( \mathfrak{D}_\infty^\M  \right)$ with $c_\M$ and $\E\left( \mathfrak{D}_\infty^\M  \right) $ which are two constants defined below.	
\end{theorem}
Our proof does not use the {\it smoothing transform} techniques, instead it relies on the {\it spine decomposition} of the branching random walk initiated by Lyons \cite{Lyo97} and a study of the tail distribution of the global minimum of the branching random walk, i.e:
$$ \M:= \inf\{ V(u),\, u\in \mathbb{T}  \}.$$
\subsection{The tail distribution of the global minimum of the branching random walk}
For $u,v\in \mathbb{T}$, we denote $u\leq v$ (resp. $u<v$) when $u$ is an ancestor of $v$ (resp. $u$ is a strict ancestor of $v$). For  $u\in \mathbb{T}-\{ \emptyset\}$, $u^-\in \mathbb{T}$ denotes the immediate ancestor of $u$, that is to say $|u^-|=|u|-1$ and $u>u^-$. Let $ \mathbb{B}(u)$ be the set of brothers of $u$, i.e
\begin{equation}
\mathbb{B}(u):=\{ v\in \mathbb{T},\, |v|=|u|\text{ and } v>u^-,\, v\neq u \}.
\end{equation}
The following theorem determines the tail distribution of the global minimum of the branching random walk. 
\begin{theorem}
	\label{MMM}
Assume (\ref{chap3criticalcondition1}), (\ref{chap3criticalcondition2}) and  (\ref{chap31.5}). 	There exists a constant $c_\M>0$ such that
	\begin{align}
	\lim_{x\to\infty} \ee^x \P\left( \M\leq -x\right) =c_\M.
	\end{align}
\end{theorem}

The link between $\M$ and $D_\infty$ can be established thanks to a {\it decompostion} of the Derivative martingale through the path associated to the vertex reaching the global minimum. Indeed for any $u\in \mathbb{T}$, by using (\ref{immedia}), one easily deduces that 
\begin{equation}
\label{theDecomp}
D_\infty= \sum_{k=1}^{|u|} \sum_{v \in \B(u_k)} e^{-V(v) }   D_\infty^{(v)} + e^{-V(u)}D_\infty^{(u)},\qquad \P-\text{ a.s.}, 
\end{equation}
with $D_\infty^{(v)}:= \lim_{n\to\infty} \sum_{u\geq v,\, |u|=n+|v|}[V(u)-V(v)] e^{-[V(u)-V(v)]}$. By the branching property of the branching random walk, one trivially checks that the random variables $(D_\infty^{(v)})_{v\in \B(u_k),k\in [1,|u|]}, D_\infty^{(u)}$ are independent copies of $D_\infty$. When $\mathtt{u}\in \mathbb{T}$ is the vertex such that $V(u)=\M$ (if several such a vertex $\u$ exist one chooses one at random among the youngest one), it gives $\P$ almost surely
\begin{align}
\label{deffrakD}  D_\infty&  = e^{-\M} \sum_{k=1}^{|\u|} \sum_{v \in \B(\u_k)} e^{\M-V(v) }   D_\infty^{(v)} + e^{\M-V(\u)}D_\infty^{(\u)}:= \ee^{-\M} \mathfrak{D}^\M
\end{align}
This time, the random variables $(D_\infty^{(v)})_{v\in \B(\u_k),k\in [1,|\u|]}, D_\infty^{(\u)}$ are independent conditionally to \\
$(V(\u_k), (V(v))_{v\in \B(\u_k)})_{ k\in [1,|\u|]}$ and distributed as the law of $D_\infty $ conditionally to $\M+V(v) \geq V(\u)$.

 \begin{theorem}
 	\label{PP}
 	Assume (\ref{chap3criticalcondition1}), (\ref{chap3criticalcondition2}) and  (\ref{chap31.5}). For any $x\in \r^+$, we denote by $(\mathfrak{D}^\M_x, \M+x)$ a random variable which has the law of $(\mathfrak{D}^\M,\M+x)$ conditionally to $\{ \M\leq -x  \}$. Then there exists a couple $(\mathfrak{D}^\M_\infty, -U)$ such that
 	\begin{align}
 	\label{eqPP}
 	\lim_{x\to\infty }(\mathfrak{D}_x^\M,\M+x) \overset{(d)}{=} (\mathfrak{D}^\M_\infty, -U).
 	\end{align}
 	Moreover $U$ is an exponential random variable with parameter $1$ independent of $\mathfrak{D}^\M_\infty$. 
 \end{theorem}
 \begin{Remark}
 	Following the proof of Theorem \ref{tailMin} below one could extract (see \eqref{defunctio}) an explicit, but rather heavy, expression for the law  $ \mathfrak{D}^\M_\infty$.
 \end{Remark}

%Finally mention that the convergence (\ref{TiTail}) may to obtain precise estimates on the modulus of continuity of the measure induced by the {\it derivative martingale}. For instance it was proven in \cite{BKNSW14} (see Theorem ??) for the "multiplicative cascade", i.e the branching random walk with independent Gaussian standard displacement. By reproducing the method of their proof, one can extend their result to every branching random satisfying (\ref{chap3criticalcondition1}), (\ref{chap3criticalcondition2}) and  (\ref{ourAssum}). In \cite{BKNSW15}, the same authors obtained similar results for a particular log-correlated model. In the last section of the present paper we shall show how to extend to a general class of log-correlated models their result by mean of the Theorem \ref{CHAOStail}.  
 
Moreover as the following proposition illustrates it, the limiting random variable $\mathfrak{D}^\M$ owns good properties of integrability.
  \begin{Proposition}
  	\label{belowplus}
  	Assume (\ref{chap3criticalcondition1}), (\ref{chap3criticalcondition2}) and  (\ref{ourAssum}), then
  	\begin{equation}
  	\label{below2plus}
  	\sup_{x\in \r^+} \E\left(     \mathfrak{D}_x^\M    \ln^2 \left(1 + \mathfrak{D}_x^\M  \right) \right)<\infty.
  	\end{equation}
  \end{Proposition}
  \begin{Remark}
  One could check that if the random variables $X$ and $\tilde{X}$ would admit finite moments of higher order, then $\mathfrak{D}^\M_x$ would also admit finite moments of higher order. 
  \end{Remark}
The proof of Theorem \ref{ourasi} is a combination of Theorem \ref{PP} and Proposition \ref{belowplus}:
  
 \noindent{\it Proof of Theorem \ref{ourasi}.} 
Recall that $ D_\infty := \ee^{-\M} \mathfrak{D}^\M$, for any $ x\geq A\geq 0$, one has
 \begin{align*}
 \P\left( D_\infty \geq \ee^x \right)&=  \P\left( \M +x\leq \ln \mathfrak{D}^\M  \right)
 \\
 & =   \P\left( \M+ x\leq  \ln \mathfrak{D}^\M ,\,  \ln \mathfrak{D}^\M \leq A \right) +  \P\left( \M +x\leq  \ln \mathfrak{D}^\M ,\,  \ln \mathfrak{D}^\M > A \right).
 \end{align*}
By decomposing the second probability on the event $\underset{p\geq 0}{\cup} \{A+p  \leq  \ln \mathfrak{D}^\M \leq  A+p+1  \}$, by Lemma \ref{tight} one has
 \begin{align*}
e^x\P\left( \M +x\leq  \ln \mathfrak{D}^\M ,\,  \ln \mathfrak{D}^\M > A \right) &\leq \sum_{p\geq 0} e^x  \P\left( \M +x\leq  A+p+1 ,\, A+p  \leq  \ln \mathfrak{D}^\M \leq  A+p+1 \right) 
\\
&\leq  \sum_{p\geq 0}  e^{A+p+1} \sup_{x\in \r} \P\left(  A+p  \leq  \ln \mathfrak{D}^\M  \Big|\M +x\leq  A+p+1\right) 
\\
&\leq  \sum_{p\geq 0} \frac{ e^{A+p+1}}{e^{A+p}(A+p)^2} \sup_{x\in \r^+} \P\left( \mathfrak{D}_x^\M \ln^2(\mathfrak{D}_x^\M)  \right) \leq c A^{-1},
 \end{align*} 
where in the last line we used Proposition \ref{belowplus}. We deduce that 
\begin{align*}
\lim_{A\to\infty} \lim_{x\to \infty} e^x\P\left( \M +x\leq  \ln \mathfrak{D}^\M ,\,  \ln \mathfrak{D}^\M > A \right)=0.
\end{align*}
On the other hand, by applying Theorem \ref{PP}, for any $p\geq 0$ we have
 \begin{align*}
 \lim_{A\to \infty }\lim_{x\to\infty } e^x\P\left(   \M+ x\leq \ln \mathfrak{D}^\M ,\,    \ln \mathfrak{D}^\M \leq A \right) &= \lim_{A\to \infty } c_\M e^{A} \P\left(  - U+ A\leq \ln \mathfrak{D}^\M_\infty ,\,  \ln \mathfrak{D}^\M_\infty \leq A  \right)
 \\
 &= \lim_{A\to\infty } c_\M  \E\left( \mathfrak{D}_\infty^\M 1_{\{ \mathfrak{D}_\infty^\M  \leq e^{A} \}} \right)
 \\
 &= c_M \E\left(\mathfrak{D}_\infty^\M  \right).
 \end{align*}
 It concludes the proof of the Theorem \ref{ourasi}.\hfill$\Box$

\paragraph{}The paper is organized as follows: Section 2 contains known facts on the {\it spine decomposition} of the branching random walk and the renewal function associated to the law of the {\it spine}.  Section 3 is devoted to the proof Theorem \ref{PP} whereas the proof of Proposition \ref{belowplus} is in section 4.

 \paragraph{Convention:} Throughout the paper, $c,\, c',\, c''$ denote generic constants and may change from paragraph to paragraph.

 \section{Preliminaries}
 \label{prelimaries}
 \subsection{Spine decomposition}
For $a\in \r$, we denote  $\P_a$ the probability distribution associated to the branching random walk starting from $a$, and $\E_a$ the corresponding expectation. Under (\ref{chap3criticalcondition1}) and (\ref{chap3criticalcondition2}), one can define the random distribution induced by
\begin{equation}
\label{defSS}
\E(f(X))=\E\left(\underset{|z|=1}{\sum}\ee^{-V(z)}f(V(z))\right),\qquad \text{for all non-negative function } f.
\end{equation}
By (\ref{chap3criticalcondition2}) we have $\sigma^2:=\E[X^2]<+\infty$. Let $(X_i)_{i\in \N^*}$ be a i.i.d sequence of copies of $X$ and for any $n\in \N$, write $S_n:=\underset{0<i\leq n}{\sum}X_i$ the mean-zero random-walk starting from the origin.
\begin{Lemma}[Biggins-Kyprianou]
	Under (\ref{chap3criticalcondition1}) and (\ref{chap3criticalcondition2}), for any $n\geq1$ and any measurable function $g:\r^n\to[0,+\infty)$,
	\begin{equation}
	\label{manytoone}
	\E\left(\underset{|z|=n}{\sum}g(V(z_1),...,V(z_n))\right)=\E\left(\ee^{S_n}g(S_1,...,S_n)\right),\qquad \text{(Many-to-one formula)}.
	\end{equation}
\end{Lemma}
Formula (\ref{manytoone}) is also a consequence of Proposition \ref{lyons} below. 

Let $\hat{L}$ be a point process which has Radon-Nikodym derivative $\int\ee^{-x}L(dx)$ with respect to the law of $L$. Conditionally to $\hat{L}= (V(z),\, |z|=1)$, let $\mathtt{w}$ be a vertex chosen among $\{z,\, |z|=1\}$ with the weights $\ee^{-V(z)}$.

 By the Kolmogorov extension Theorem there exists a probability measure $\Q$ such that for any $n\geq 0$,
\begin{equation}
\Q|_{\mathcal{F}_n}:=W_n\bullet\P|_{\mathcal{F}_n},
\end{equation}
where $\mathcal{F}_n$ denotes the sigma-algebra generated by the positions $(V(z),\,|z|\leq n)$ up to time $n$.  Lyons \cite{Lyo97} gave the following description of the branching random walk under $\Q$: 
\begin{enumerate}
\item  First consider $(\mathtt{w}_i,\xi_{i},\hat{L}_i)_{i\geq 1}$ an i.i.d sequence with $(\mathtt{w}_1,\xi_{i},\hat{L}_1)\overset{(d)}{=} (\mathtt{w},V(\mathtt{w}),\hat{L})$.

\item Then for any $i\in \N^*$, attach to $\mathtt{w}_i$ the point process $\hat{L}_i$, in order to have
\begin{enumerate}
	\item an infinite {\it spine}, denoted by  $(w_n)_{n\in \N}$ ($w_0=\emptyset$, $w_1=\emptyset\mathtt{w}_1$, $w_2=\emptyset \mathtt{w}_1\mathtt{w}_2$ ...) 
	
	\item and the set of immediate brothers of the vertices of the spine:  $(\B(w_i)=\{ u>w_{i-1},\, |u|=i,\, u\neq w_i  \})_{i\in \N^*}$.  
\end{enumerate}

\item Finally for any $i\in \N^*$ and $u\in \B(w_i)$ attach an independent Branching random walk, sampled under $\P$, rooted at $u$ and denoted by $\text{BRW}(u)$. 

\end{enumerate}
We still call $\mathbb{T}$ the genealogical tree of the process. 
\begin{Proposition}
	\label{lyons}
	Suppose (\ref{chap3criticalcondition1}) and (\ref{chap3criticalcondition2}). For any $|z|=n$, we have
	\begin{equation}
	\label{no(iii)}
	\Q\left\{w_n=z\Big| \mathcal{F}_n\right\}=\frac{\ee^{-V(z)}}{W_n};
	\end{equation}
Moreover  the spine process $(V(w_n),n\geq 0)$ has the distribution of the centered random walk $(S_n)_{n\geq 0}$ satisfying (\ref{manytoone}).
\end{Proposition}

\begin{Remark}
	The change of probability is now a standard technique. We refer to \cite{LPP95} for the case of the Galton-Watson tree, to \cite{CRo88} for the branching Brownian motion, and to \cite{BKy04} for the spine decomposition in various types of branching.
\end{Remark}

\paragraph{A time reversal identity:} Under $\Q$ the branching random walk is constructed uniquely thanks to the i.i.d sequence $(\mathtt{w}_i,\xi_{i},\hat{L}_i)_{i\geq 1}$ and the independent branching random walk attached to each brothers of the spine. In particular the vectors $(\mathtt{w}_i,\xi_{i},\hat{L}_i)_{i\in [|1,k|]}$ and $(\mathtt{w}_{k-i},\xi_{k-i},\hat{L}_{k-i})_{i\in [|1,k|]}$ have the same law. It induces the following {\it time reversal} identity
\begin{align}
\nonumber &\E_\Q\left(  \varphi\left[   (V(w_i),  (V(u),BRW(u))_{   u\in \B(w_i)}  )_{i\in [|1,k|]}\right]   \right) 
\\
\label{reversal}&\qquad\qquad =\E_\Q\left(  \varphi  \left[    ( V(w_k)-V(w_{k-i }) ,  (V(u), BRW(u))_{   u\in \B(w_{k-i +1}) }     )_{i\in [|1,k|]} \right]  \right).
\end{align}
which holds for any continuous and bounded functional $\varphi$. This identity will be crucial in the end of the proof of Theorem \ref{tailMin} (Subsection \ref{proffttailmin}).

 \paragraph{The probability $\Q_k\otimes\P$:} Finally, for any $k\in \N$ we introduce the probability $\Q_k\otimes\P$ under which the branching random walk up to time $k$ is distributed as a branching random walk under $\Q$ and after the time $k$ every alive particle at time $k$ will branch according to the original point process $L$ under $\P$.

 \subsection{The renewal function associated to a one-dimensional random walk}
Thanks to the spine decomposition technique many questions concerning the whole branching random walk, can be reduced in one computation involving the standard random walk $(S_n)_{n\geq 0}$ introduced in (\ref{defSS}).  We collect here some known facts on the renewal function and the paths of such a standard random walk.

Recall that $\E(S_1)=0$ and $\sigma^2=\E( S_1^2)\in(0,\infty)$.  Let $(H_i^-)_{i\geq 0}$ and $(H_i^{+})_{i\geq 0}$ are respectively the strict descending and ascending ladder height of $ (S_n)_{n\geq 0}$. It means that $H_i^{-}=S_{T_i^-}$ and $H_i^+=S_{T_i^+}$ with $T_0^-=T_0^+=0$ and $T_{i+1}^-:= \inf\{ j>T_i^-,\, S_j<S_{T_i^-}\}$ and $T_{i+1}^+ := \inf\{ j>T_i^+,\, S_j>S_{T_i^+}\}$. According to Feller \cite{Fel71},  $\E(|H_1^+|)<+\infty$ and $\E(|H_1^-|)<+\infty  $, so we can define the renewal functions associated to $ (S_n)_{n\geq 0}$ by
 \begin{equation}
 	\label{definih0}
 	R^-(u):=\sum_{j\geq 0} \P\left( H_j^-\geq -u\right),\qquad R^+(u):= \sum_{j\geq 0} \P\left( H_j^+\leq u\right),\quad u\geq 0
 \end{equation}
By using the time reversal property of $ (S_n)_{n\geq 0}$ we can rewrite these two functions as 
 \begin{equation}
 	R^-(u)=\underset{j\geq 0}{\sum}\P\left( \underset{i\in [|1,j|]}{\max}\, S_i< 0,\, S_j\geq -u  \right),\qquad   	R^+(u)=\underset{j\geq 0}{\sum} \P\left( \underset{i\in [|1,j|]}{\min}\, S_i> 0,\, S_j\leq u  \right)   ,\qquad  u\geq 0,
 \end{equation}
 %which is also equal by the duality lemma to
 %\begin{equation}
 %\label{chap32.3}
 %h_0(u):=\overset{\infty}{\underset{j=0}{\sum}}\P\left( S_1< 0,\,...,S_j< 0,\, S_j\geq -u\right),\qquad u\geq 0,
 %\end{equation}
with the conventions: $\max_{i\in \emptyset} S_i=-\infty$, $\min_{i\in \emptyset} S_i=+\infty $. Observe that $R^+$ and $R^-$ are increasing and $R^+(0)=R^-(0)=1$, moreover for any $u\geq 0$, $R^+$ and $R^-$ satisfy
 \begin{equation}
 	\label{chap33.40}
 	R^-(u)=\E\left(R^-(S_1+u)\1_{\{S_1\geq -u\}}\right),\qquad R^+(u)=\E\left(R^+(u-S_1)\1_{\{S_1\leq u\}}\right).
 \end{equation}
According to the Theorem 1, Section XVIII.5 p.612 in \cite{Fel71}, there exists $C^-,\, C^+>0$ such that 
 \begin{equation}
 	\label{chap3lldda}
 	C^-:=\underset{u\to\infty}{\lim}\frac{R^-(u)}{u},\qquad  	C^+:=\underset{u\to\infty}{\lim}\frac{R^+(u)}{u},\qquad \forall u\geq 0,
 \end{equation}
and furthermore by the Blackwell renewal theorem (see for instance Theorem 4.4.3 in \cite{Dur10}), for any $h>0$, 
\begin{align}
\label{blackwi}
\lim_{u\to\infty } R^-(u+h)-R^{-}(u) =C^-h,\qquad \lim_{u\to\infty} R^{+}(u+h) -R^+(u)= C^+h.
\end{align} 
As a consequence there exist constants $c_1,\,C_1>0$ such that  
 \begin{equation}
 	\label{chap32.7}
 	c_1(1+u)\leq R^-(u),\, R^+(u)\leq C_1(1+u),\qquad u\geq 0.
 \end{equation}
By Kozlov Formula (12) in \cite{Koz76}, we know also that when $n\to \infty$ uniformly in $u\in [0, (\log n)^{30}]$,
 \begin{equation}
 	\label{chap3lldda2}
 	\P\left(\min_{j\in [|0,n|]}S_j\geq -u\right)= \frac{\theta_-  R^-(u)+ o(1)}{n^\frac{1}{2}},\qquad \P\left(\max_{j\in [|0,n|]}S_j\leq u\right)= \frac{\theta_+ R^+(u)+ o(1)}{n^\frac{1}{2}},
 \end{equation}

Mention also an inequality due to \cite{AShi10}{\string:} {\it there exists $c>0$ such that for $u>0,a\geq0,\,b\geq 0$ and $n\geq 1$,
   	\begin{equation}
   		\label{chap3Aine}
   		\P\left(\min_{j\leq n}S_j\geq -a,\,b-a\leq S_n\leq b-a+u\right)\leq c\frac{(u+1)(a+1)(b+u+1)}{n^{\frac{3}{2}}}.
   	\end{equation}}
 Finally we recall one useful result proved in \cite{Aid11}
 \begin{Lemma}[\cite{Aid11}]
 	\label{techniqAID}
 	Let $a>0$. There exists a constant $c(a)>0$ such that for any $z\geq 0$,
 	\begin{align}
 	\E_z\left[ \sum_{l\geq 0}\ee^{-a S_l} 1_{\{ \min_{j\leq l} S_j\geq 0 \}}   \right]= c(a)<\infty.
 	\end{align}
 \end{Lemma}
 
\subsection{The renewal function starting from any point}
Let us introduce the following extension of $R^{-}$ which will be ubiquitous through the paper:
 \begin{equation}
 \label{tildiR}
 \tilde{R}(x,a):= \sum_{j\geq 0} \P_{-a}\left( \max_{i\in [|1,j|]}S_i< 0,\, S_j\geq -x   \right),\quad \forall a,x\geq 0.
 \end{equation}
 Remark that $\tilde{R}(x,0)= R^-(x)$. 
 \begin{Lemma}
 	\label{renouvDoux}
 	For any $x, a>0$ we have
 	\begin{equation}
 	\label{eqrenouvDoux}
	\tilde{R}(x,a)  =  \theta_0 R^-(x)\left\{  R^+(a) -K_a  \right\}  +   \theta_0 \int_{x-a}^{x}  \left\{K_{a-x+u}-R^+(a-x+u) \right\}  dR^-( u).
 	\end{equation}
 	where $R^-$ and $R^+$ are the renewal functions defined in (\ref{definih0}) and for any $u\in \r$,  $K_u:= \E\left(\sum_{i\geq 0} 1_{\{ H_i^+=u \}}\right)$  and $\theta_0:= \sum_{j\geq 0} \P\left( \max_{l\in [|0,j|]}S_l\leq 0,\, S_j=0 \right)$.
  \end{Lemma}
 We stress that the formula \eqref{eqrenouvDoux} is not true for $a=0$. Furthermore remark that by the Blackwell renewal theorem (\ref{blackwi}), we also have
\begin{equation}
\label{laRemark}
	 	\lim_{x\to\infty}  \theta_0 \int_{x-a}^{x}  \left\{K_{a-x+u}-R^+(a-x+u) \right\}  dR^-(u)= C^-\int_0^a (K_u- R^+(u))du.
\end{equation}
 In particular we shall use this Lemma in combination with (\ref{laRemark}) at the end of the proof of Theorem \ref{tailMin}.

 {\it Proof of Lemma \ref{renouvDoux}.} According to the time reversal property of the random walk $(S_i)_{i\in [|0,j|]}$, for any $a>0$ we have 
 \begin{align*}
 \tilde{R}(x,a)&= \sum_{j\geq 0} \P\left( S_j< a+ \min_{i\in[|0,j-1|]}S_i,\, S_j\geq a-x    \right)
 \\
 &=  \E\left(\sum_{j\geq 0} 1_{\{ H_j^-\geq -x  \}} \times     \sum_{i\geq 0} 1_{\{ a-x-H_{j}^-  \leq S_{T^-_j +i}-H^-_j < a,\, \underset{l\in [|0,i|]}{\min} S_{T^-_j+l}\geq  H^-_j   \}}  \right)  .
 \end{align*}
 Observe that for every $j\in \N$ the sequence of paths $(S_{T^-_j+i}-H^-_j)_{i\in [|0,T_{j+1}^- -T_j^-|]}$ are independent and identically distributed as an excursion above $0$ stopped when it reaches $(-\infty,0)$. Then we have 
\begin{align*}
 \tilde{R}(x,a)&= \sum_{j\geq 0} \E\left( \1_{\{H^-_j\geq -x\}}   \times \sum_{i\geq 0}\P( a-x-H_j^-\leq S_i < a,\, \min_{l\in [|0,i|]}S_l\geq 0) \right).
\end{align*}
By the time reversal property of $(S_n)_{n\geq 0}$, 
 \begin{align*}
 \sum_{i\geq 0}\P\left( S_i < a,\, \underset{l\in [|0,i|]}{\min}S_l\geq 0  \right)&= \sum_{i\geq 0} \P\left(  S_i<a,\, S_i\geq \underset{l\in [|0,i|]}{\max} S_l\right)
 \\
 &= \E\left(  \sum_{i\geq 0} \1_{\{ H_i^+<a \}} \sum_{j\geq 0} \1_{\{  \underset{l\in [|0,j|]}{\max} S_{T_i^++l}-H_i^+\leq 0,\, S_{T_i^++j}-H_i^+=0 \}}  \right)
 \\
 &=\E\left(\sum_{i\geq 0} \1_{\{ H_i^+<a \}} \right) \times \sum_{j\geq 0} \P\left( \max_{l\in [|0,j|]}S_l\leq 0,\, S_j=0 \right)
 \\
 &= \{ R^+(a) - K_a  \} \theta_0 .
 \end{align*}
 with $K_a:= \E\left( \sum_{ i\geq 0} 1_{\{ H_i^+=a\}}\right)$ and $\theta_0= \sum_{j\geq 0}\P\left(  \max_{l\in [|0,j|]}S_l\leq 0,\, S_j=0 \right)$. Thus we deduce that
\begin{align*}
 \tilde{R}(x,a)&= \theta_0\sum_{j\geq 0} \E\left( \1_{\{H^-_j\geq -x\}}   \times \left\{ R^{+}(a)- R^{+}(a-x-H_j^-) +K_{a-x-H_j^-}- K_a   \right\}  \right) 
 \\
 &=\theta_0 R^-(x)\left\{  R^+(a) -K_a  \right\}  + \theta_0 \int_{x-a}^{x}  \left\{K_{a-x+u}-R^+(a-x+u) \right\}  dR^-(u) .
% &= R^-(x)\{ R^+(a) - \sum_{i\geq 0} \P\left( H_i^+=a  \right) \}\sum_{i\geq 0}\P\left(  \max_{l\in [|0,i|]} S_l\leq 0,\, S_i=0\right).
\end{align*}
 It concludes the proof of the Lemma \ref{renouvDoux}. \hfill$\Box$
 \\

We end this section by the following useful bound on $\tilde{R}(x,a)$:
 \begin{Lemma}
 	\label{Lemconsequec}
There exists $c>0$ such that for any $x,a,b\geq 0$,
 	 \begin{align}
 	 \label{consequencese}
 	 \tilde{R}(x+b,a)-\tilde{R}(x,a)\leq c(1+ a)(1+b)^2.
 	 \end{align}
 \end{Lemma}
\noindent {\it Proof of Lemma \ref{Lemconsequec}.} Let $\tau_{a-x}:= \inf\{ k\geq 0,\, S_k< a-x \}$ be a stopping time. By the definition of $\tilde{R}$ in (\ref{tildiR}), 
\begin{align*}
\tilde{R} (x+b,a) -\tilde{R} (x,a)&= \sum_{j\geq 0}   \P\left( S_j< a+ \min_{i\in[|0,j-1|]}S_i,\,  a-x > S_j\geq a-(x+b)   \right)
\\
&= \E\left( \sum_{j\geq \tau_{a-x}} 1_{\{  S_j< a+ \min_{i\in[|0,j-1|]}S_i,\,  a-x > S_j\geq a-(x+b)     \}} \right) . 
\end{align*}
By the Markov property at time $\tau_x$ on has
\begin{align*}
\tilde{R} (x+b,a) -\tilde{R} (x,a) &\leq \E\left(  \sum_{j \geq 0}   \P\left( S_j< a+ \min_{i\in[|0,j-1|]}S_i,\,  a-x > S_j +z\geq a-(x+b)        \right)_{\big|  z=S_{\tau_{a-x}} }  \right)
\\
&\leq \sup_{-x-b \leq   z \leq a-x} \sum_{j\geq 0} \P\left(    S_j< a+ \min_{i\in[|0,j-1|]}S_i,\,  a-x > S_j +z\geq a-(x+b)         \right)
\\
&=\sup_{  b+ a\geq  z  \geq 0}   \sum_{j\geq 0} \P\left(    \max_{i\in [|1,k|]}S_i< a,\,  z> S_j  \geq z-b         \right) \leq   c(1+a) (1+b)^2.
\end{align*}
where in the last line we operated a time reversal then used (\ref{chap3Aine}).
\hfill$\Box$

\section{The derivative martingale seen from the global minimum}

For any $j\in \N^*$, let us denote $\sum_{u\in \B(w_j)} \delta_{   \zeta_u^{(j)}}:=\sum_{u>w_{j-1},\,u\neq w_j,\,  |u|=j}\delta_{{\{ V(u)-V(w_{j-1}) \}}}$ the point process formed by the position the brothers of $w_j$. We introduce the truncated version of $\mathfrak{D}^\M$, i.e
\begin{align}
\label{UPZE}
\forall u\in \mathbb{T},\qquad \mathfrak{D}^{u,\geq t}:=e^{V(u)}\sum_{j=|u|-t}^{|u|} e^{-V(u_{j-1})} \sum_{v \in \B(u_j)} e^{- \zeta_v^j}   D_\infty^{(v)} + D_\infty^{(u)}.
\end{align}
%In (\ref{trivago}) we defined $\Q_r$, the law of the point process $\mu(w_1)$ conditioned to have $V(w_1)=r$. For any $\theta \in \r^+$, $A\subset \r^+$ compact set, let us define 
%\begin{align}
%\label{DEFPSI}
%\psi^{(\theta,A)}(r,y,z)&:= \E_{\Q_r}\left(  \prod_{u\in \mathbb{B}(w_1)} \E_\P\left(1_{\{\mathtt{M} \geq z-V(u)    \}} \ee^{-\theta \eta^{\mathtt{M}}(A,V(u)+y, z-V(u))}   \right)\right)
%\\
%\nonumber &= \E_{\Q_r\otimes \P}\left(  \prod_{u\in \mathbb{B}(w_1)}  1_{\{\mathtt{M}^{(u)} \geq z-V(u)    \}} \ee^{-\theta \eta^{(u),\mathtt{M}}(A,V(u)+y, z-V(u))}    \right) ,\qquad r,\, y,\, r\in \r^3.
%\end{align}
%with $\eta^{(u),\mathtt{M}}(\cdot,\cdot,\cdot))_{u\in \mathbb{T}}$ an i.i.d family of random variable, distributed as
%\begin{align}
%\eta^{\mathtt{M}}(A,x,y):= \sum_{v\in \mathbb{T}}\1_{\{  \underset{j\in [|0,|v||]}{\min }  {V}(v_j)= V(v),\, x+V(v)\leq 0,\, V(v)-y\in A \}}
%\end{align} 
	\begin{theorem}
		\label{tailMin}
Assume (\ref{chap3criticalcondition1}), (\ref{chap3criticalcondition2}) and  (\ref{chap31.5}). 	Let $t \in \N^*$ be an integer. Let $\mathtt{u}\in \mathbb{T}$ be the vertex such that $V(u)=\M$ (if several such a vertex $\u$ exist one chooses one at random among the youngest one). There exists a non-null functional $\mathcal{E}_t$, such that for any continuous and bounded function $\varphi: \r\mapsto \r^+$ we have the following limit 
\begin{align*}
\lim_{x\to \infty} \ee^{x} \E\left(  \varphi\Big(\mathfrak{D}^{\u,\geq t}   \Big)\1_{\{  \M\leq -x \}} \right):=  \mathcal{E}_t(\varphi).
\end{align*}

	\end{theorem}
\begin{Remark}
	\begin{enumerate}
	\item By taking $\varphi$ constant equal to $1$, it implies Theorem \ref{MMM}.
	
	\item An explicit expression of the functional $\mathcal{E}$ is written in (\ref{defunctio}).
\end{enumerate}
\end{Remark}
 
\noindent{\it Proof of Theorem \ref{PP}.} By Theorem \ref{tailMin} we can affirm that for any continuous and bounded function $\varphi: \r\mapsto \r^+$ and any $t\in \N^*$, 
\begin{align}
	\nonumber &\lim_{x\to \infty}  \E\left(  \varphi\Big(\mathfrak{D}^{\u,\geq t}\Big)  \1_{{\{ \M+x \leq -u  \}}}  \big|  \M\leq x\right)=  \frac{ e^{-u} \mathcal{E}_t(\varphi)}{ \mathcal{E}_t(1)}= \frac{ e^{-u} \mathcal{E}_t(\varphi)}{ \mathcal{E}_0(1)}.
\end{align}
Moreover by Proposition \ref{belowplus} the family of distribution of $  ( \mathfrak{D}^{\u,\geq t}   ,\M+x) $ conditionally to $\{\M\leq x   \} $ is clearly tight. By applying the classical Lévy' Theorem, there exists a couple of independent random variables $(\mathfrak{D}^{\u,\geq t}_\infty , U)$ such that conditionally to $\M\leq -x$, 
\begin{align}
( \mathfrak{D}^{\u,\geq t}   ,\M+x) \overset{\text{weakly}}{\Longrightarrow } (\mathfrak{D}^{\u,\geq t}_\infty  , -U),
\end{align}
with $U$ an exponential random variable with parameter $1$. We now are in shape to prove the convergence (\ref{eqPP}). Indeed as the family of distribution $(\mathfrak{D}^\M_\infty,\M+x)$ conditionally to $\{\M\leq -x\}$ is tight, it suffices to prove that for any $\theta_1,\theta_2\in \r^+$,
\begin{align}
\label{termi}\lim_{x\to\infty} \E\left(  e^{-\theta_1 \mathfrak{D}^\M + \theta_2 (\M+x) } \Big| \M\leq -x \right)= \lim_{t\to +\infty} \frac{\mathcal{E}_t(e^{-\theta_1 \cdot })}{\mathcal{E}_0(1)} \frac{\theta_2}{1+\theta_2}.
\end{align}
Notice that the right hand limit term exists as $t\mapsto \mathcal{E}_t(\ee^{-\theta_1 \cdot})$ is decreasing and positive. By Lemma \ref{toreste} we clearly have, for any $\epsilon >0$
\begin{align*}
\lim_{t\to\infty} \lim_{x\to\infty}\E\left( e^{-\theta_1 \mathfrak{D}^{\u,\geq t} + \theta_2 (\M+x) }- e^{-\theta_1 \mathfrak{D}^\M + \theta_2 (\M+x) }   \Big| \M\leq -x \right)\leq \theta_1\epsilon + \lim_{t\to\infty} \lim_{x\to\infty} \P\left(   \mathfrak{D}^\M  - \mathfrak{D}^{\u,\geq t} \geq \epsilon   \Big| \M\leq -x \right)\leq \theta_1 \epsilon,
\end{align*}
which suffices to obtain (\ref{termi}) and concludes the proof of  Theorem \ref{PP}.

\begin{Remark} Following step by step the proof of Theorem \ref{tailMin} it is plain to check the existence of a non-null functional $\tilde{\mathcal{E}}_t$ such that for any continuous and bounded function $\varphi: \r^t\mapsto \r^+$,
	\begin{align*}
	&\lim_{x\to \infty} \ee^{x} \E\left(  \varphi\Big(  V(|\u|_{|\u|-1})-\M,..., V(\u_{|\u|-t})-\M   \Big) 1_{\{  \M\leq -x    \}   }\right) =   \tilde{\mathcal{E}}_t(\varphi) .
	%			\\
	%			\nonumber	& C^-\theta_0     \E_{\bar \Q}\left(  \1_{\{ \underset{j\in [|1,+\infty|]}{\max}V(w_j)<0\}} \frac{  \prod_{j\geq 0} \prod_{u\in \B(w_{j+1})} \1_{\{\M^{(u)}\geq \Xi^{(j+1)}_u   ,\, \M_{j-1}^{(u)} >   \Xi^{(j+1)}_u \}}      \varphi\Big(  -V(w_{1}),..., -V(w_{t})     \Big) }{1+    \sum_{j=0}^{\infty} \sum_{u\in {\B}(w_{j+1})}  N_j^{(u)}(\Xi^{(j+1)}_u )  }         \right) .
	\end{align*}
	It would prove the convergence, when $x\to\infty$ of the distribution of  $(V(\u_{|\u|-1})-\M,..., V(\u_{|\u|-t})-\M )$ conditionally to $\M\leq -x$. 
	
\end{Remark}

\subsection{Upper and lower bound for the tail distribution of $\mathtt{M}$}
The following Lemma ensures that the constants $\mathcal{E}_0(1)$ of Theorem \ref{tailMin} is non null. 
\begin{Lemma}
\label{tight}
Assume (\ref{chap3criticalcondition1}), (\ref{chap3criticalcondition2}) and  (\ref{chap31.5}). There exists $c_1>0$ such that for any $x> 0$, 
\begin{equation}
\label{eqticht}
c_1\leq \ee^{x} \P\left( \M\leq -x \right)\leq 1.
\end{equation}
\end{Lemma}
\noindent{\it Proof of Lemma \ref{tight}.} We recall here the proof of the upper bound written in \cite{Aid11},
\begin{eqnarray*}
\P\left( \mathtt{M}\leq -x \right)& \leq& \E\left( \sum_{k\geq 1} \sum_{|u|=k} 1_{\{  \min_{i\leq k-1} V(u_i) >-x,\,  V(u)  \leq -x  \}} \right)
\\
&=& \sum_{k\geq 1} \E( \ee^{S_k} 1_{\{ \min_{i\leq k-1}S_i>-x,\, S_k\leq -x \}})
\\
&\leq &\sum_{k\geq 1}\P(  \min_{i\leq k-1}S_i>-x,\, S_k\leq -x )\ee^{-x}\leq \ee^{-x}.
\end{eqnarray*}
To prove the lower bound we will use the second moment method and the idea of {\it good vertex} first introduced by Aïdékon in \cite{Aid11}. It consists to exclude the vertices of the branching random walk which make explode the second moment. For any $x,L>0$, let us define
$$ N_L(x):=  \sum_{k\geq 1} \sum_{|u|=k} 1_{\{  \min_{i\leq k-1} V(u_i) > V(u),\, V(u) \in I(x),\, \sum_{j=1}^k\sum_{u\in \B(u_j)}\ee^{-V(u)-x}\leq L \}}, $$ 
with 
\begin{align*}
I(x):=    [-x-1,x) .
\end{align*}
By the Paley-Zygmund inequality, note that for any $x,L> 0$,
$$\P(\M\leq  -x)\geq \P\left( N_L(x)>0 \right)\geq \frac{\E( N_L(x))^2}{\E(N_L(x)^2)}.$$
To prove a lower bound on $\E\left( N_L(x)\right)$, observe that
\begin{align}
\E(N_L(x))= \E\left(N_\infty(x) -\bar{N}_L(x)\right),
\end{align}
with 
$$ \bar{N}_L(x):=  \sum_{k\geq 1} \sum_{|u|=k} 1_{\{  \min_{0\leq i\leq k-1} V(u_i) > V(u),\, V(u) \in I(x),\, \sum_{j=1}^k\sum_{u\in \B(u_j)}\ee^{-V(u)-x}>L \}}. $$
Moreover by the Proposition \ref{lyons}
\begin{align}
\nonumber \E\left(  N_\infty(x) \right) = \E\left( \sum_{k\geq 1} \ee^{S_k} \1_{\{ \min_{0\leq i\leq k-1} S_i> S_k,\, S_k\in I(x)  \}}\right) &\geq \ee^{-x-1} \sum_{k\geq 1}  \P\left(   \min_{0\leq i\leq k-1} S_i> S_k,\, S_k\in I(x)    \right)
\\
\label{toutou0}&= [ R^-(x+1)-R^-(x)  ]\ee^{-x}\geq c\ee^{-x},
\end{align}
where we used \eqref{blackwi} in the last inequality. On the other hand, again by the Proposition \ref{lyons}
\begin{align*}
\E\left( \bar{N}_L(x)\right)& = \sum_{k\geq 1} \E\left( \sum_{|u|=k} 1_{\{  \min_{0\leq i\leq k-1} V(u_i) > V(u),\, V(u) \in I(x),\, \sum_{j=1}^k\sum_{u\in \B(u_j)}\ee^{-V(u)-x}>L \}}  \right)
\\
&\leq \ee^{-x} \sum_{k\geq 1} \Q\left( \min_{0\leq i\leq k-1} V(w_i) > V(w_k),\, V(w_k) \in I(x),\, \sum_{j=1}^k\sum_{u\in \B(w_j)}\ee^{-V(u)-x}>L \right) 
\\
&\leq \ee^{-x} \sum_{k\geq 1} \sum_{j=1}^k  \Q\left( \min_{0\leq i\leq k-1} V(w_i) > V(w_k),\, V(w_k) \in I(x),\,  \Delta_j e^{-V(w_{j-1})} >\frac{L}{c_\kappa {(k-j+1)}^{\kappa}} \right) 
\\
&:= \ee^{-x} \sum_{k\geq 1} \sum_{j=1}^k  (1)_{j,k}.
\end{align*}
with $\Delta_j= \sum_{u\in \B(w_j)}\ee^{-[V(u)-V(w_{j-1})]}$, $\kappa>1$ and $c_{\kappa}>0$ large enough to have $\sum_{j\geq 1}j^{-\kappa}\leq c_\kappa$. The random variables $(V(w_j)-V(w_{j-1}), \Delta_j)_{j\in [|1,k|]}$ are independent and identically distributed, thus 
\begin{align*}
\left(V(w_j)-V(w_{j-1}), \Delta_j\right)_{j\in [|1,k|]} \overset{\text{(d)}}{=} \left( V(w_{k-j+1})-V(w_{k-j}), \Delta_{k-j+1}   \right)_{j\in [|1,k|]}. 
\end{align*}
By operating a time reversal one gets
\begin{align*}
(1)_{j,k}= \Q\left( \max_{i\in [|1,k|]} V(w_i)<0,\, V(w_k)\in I(x),\, \Delta_{k-j+1} e^{V(w_{k-j+1})} >L \frac{x+e^{V(w_k)}}{c_\kappa (k-j+1))^{\kappa}}    \right) .
\end{align*} 
In others words, for any $j\in [|1,k|]$, 
\begin{align*}
 (1)_{k-j+1,k} & \leq  \Q\left( \max_{i\in [|1,k|]} V(w_i)<0,\, V(w_k) \in I(x),\, \ee^{V(w_{j})}\Delta_j>L  \frac{\ee^{V(w_k)+x}}{c_\kappa j^\kappa} \right) 
\\
&\leq \Q\left( \max_{i\in [|1,k|]} V(w_i)<0,\, V(w_k) \in I(x),\,   -V(w_{j}) +\ln \frac{L}{4c_\kappa} <  \ln \Delta_j +\kappa \ln j \right),
\end{align*}
By the Markov property at times $j$ we get
\begin{align*}
  (1)_{k-j+1,k}  &\leq   \E_\Q\left(\1_{\{  \max_{i\in [|1,j|]} V(w_i)<0,\,  -V(w_{j}) +\ln \frac{L}{4c_\kappa} <  \ln \Delta_j +\kappa \ln j   \}}  \P_{V(w_j)}\left(  \max_{i\in [|1,k-j|]}S_i<0,\, S_{k-j}\in I(x)  \right) \right).
  \end{align*}
 Then by recalling the definition (\ref{tildiR}) and reversing the indices one can affirm that
\begin{align*}
\E(\bar{N}_L(x))&\leq \ee^{-x}\sum_{k\geq 1} \sum_{j=1}^k (1)_{k-j+1,k}
\\
&= \ee^{-x} \sum_{j\geq 1}  \E_\Q\left(\1_{\{  \underset{j\in [|1,j|]}{\max} V(w_i)<0,\,  -V(w_{j}) +\ln \frac{L}{4c_\kappa} <  \ln \Delta_j  +\kappa \ln j \}} [\tilde{R}(x+1,-V(w_j)) - \tilde{R}(x,-V(w_j))]  \right)
\\
&\leq c\ee^{-x} \sum_{j\geq 1}  \E_\Q\left(\1_{\{  \max_{j\in [|1,j|]} V(w_i)<0,\,  -V(w_{j})  <  \ln \Delta_j +\kappa \ln j -\ln \frac{L}{4c_\kappa}  \}}(1 -V(w_j) ) \right).
\end{align*}
where we used Lemma \ref{Lemconsequec} in the last inequality. When $L> 4c_\kappa$, one has $ -V(w_{j})  <   \ln \Delta_j +\kappa \ln j  . $ Moreover by setting $\Delta^+_j:= \sum_{u\geq w_{j-1},\, |u|=j} \ee^{- [V(u)-V(w_{j-1})]}$ one also has
\begin{align*}
-V(w_{j-1})= -V(w_j) + V(w_j) - V(w_{j-1}) &\leq \ln \Delta_j +\kappa \ln j -\ln \frac{L}{4c_\kappa}  - \ln e^{- [V(w_j)-V(w_{j-1})]}
\\
&\leq  \ln \Delta^+_j +\kappa\ln j -\ln \frac{L}{4c_\kappa}.
\end{align*}
Let $\Delta^+$ be a random variable distributed at $\Delta^+_j$ independent of everything, by the Markov property at time $j-1$ and (\ref{chap32.7}) we get that 
\begin{align}
\nonumber\ee^x \E(\bar{N}_L(x)) & \leq c     \E_\Q\left(\sum_{j\geq 1} \Q\left( \max_{j\in [|1,j|]} V(w_i)<0,\,  -V(w_{j-1})  <   \ln \Delta^+ +\kappa \ln j -\ln \frac{L}{4c_\kappa}  \right)(1+ \ln_+ \Delta^+ + \ln j)  \right)  
\\
\nonumber &\leq   c'  \sum_{j\geq 1} \frac{(1+2\ln j  - \ln \frac{L}{4c_\kappa})_+^3}{j^{\frac{3}{2}}}  +  c' \sum_{j\geq 1}   \frac{(1+ \ln j)}{j^{\frac{3}{2}}} \E_\Q\left((1+ 2\ln_+\Delta^+- \ln \frac{L}{4c_\kappa})^2 \1_{\{  \Delta^+>    \sqrt{\frac{L}{4c_\kappa}} \}} \right)  
\\
\nonumber &\qquad\qquad\qquad\qquad  + c'\E_\Q\left(\sum_{j\geq 1} \Q\left( \max_{j\in [|1,j|]} V(w_i)<0,\,  -V(w_{j-1})  <   2\ln \Delta^+  -\ln \frac{L}{4c_\kappa}  \right)(1+ \ln_+ \Delta^+)  \right) 
\\
  &\leq c'' \E_{\Q} \left( [1+ \ln_+\Delta^+ ]^2\1_{\{  \Delta^+>    \sqrt{\frac{L}{4c_\kappa}} \}} \right)+ c''\sum_{j\geq (\frac{L}{4c_\kappa})^{1\over \kappa}}  \frac{(1+2\kappa \ln j)^3}{j^\frac{3}{2}}.\nonumber 
\end{align}
The random variable $\Delta^+$ is stochastically dominated by $X$, thus by (\ref{chap31.5}) we deduce that uniformly in $x>0$, 
\begin{align}
\label{toutou}\lim_{L\to \infty} e^{x} \E\left( \bar{N}_L(x) \right) = 0.
\end{align}
 By combining (\ref{toutou0}) and (\ref{toutou}) we deduce that for a large enough $L>0$, there exists $c>0$ such that for any $x\geq 0$, 
 $$ \E\left( N_L(x) \right) \geq c\ee^{-x}.$$
Now we shall study the second moment of $N_L(x)$. By definition 
\begin{align*}
N_L(x)^2= N_L(x) +    \sum_{n\geq 1} \sum_{|u|=n} \left( 1_{\{  ... \}}   \sum_{p\geq 1} \sum_{|v|=p}  1_{\{ ...\}} 1_{\{ v\neq u \}} \right).
\end{align*}According to the Proposition \ref{lyons}, it leads to
%\begin{align*}
% A(x)^2 &= \sum_{k_1\geq 0} \sum_{k_2\geq 0} \sum_{|u|=k_1}\sum_{|v|=k_2} 1_{\{  \min_{j\leq k_1-1} V(u_j)>V(u),\, V(u)\in I(x)   \}}1_{\{  \min_{j\leq k_2-1} V(v_j)>V(v),\, V(v)\in I(x)   \}}
% \\
% &=A(x)+\sum_{k\geq 0} \sum_{|u|=k}  \left( \sum_{v \geq u^{(1)}}   1_{\{  \min_{j\leq k_1-1} V(v_j)>V(v),\, V(v)\in I(x)   \}}        \right)  \left( \sum_{w\geq u^{(2)}}   1_{\{  \min_{j\leq k_2-1} V(w_j)>V(w),\, V(w)\in I(x)   \}}\right).
%\end{align*}
%Taking the expectation we obtain
%\begin{align*}
%\E(A(x)^2) & \leq \E(A(x))+ \sum_{k\geq 0} \E\left( \sum_{|u|=k}1_{\{ \min_{j\leq k}V(u_j)>-(x+1) \}} \sum_{v\neq w,\, v\geq u,\, |v|=|w|=k+1} \ee^{-x-V(v)} \ee^{-x-V(w)} \right)
%\end{align*}
\begin{align*}
&\E\left( N_L(x)^2\right) \leq \E(N_L(x))+  \sum_{n\geq 0} \E_{\Q_n\otimes \P}\left( \ee^{V(w_n)}     \1_{\{  \underset{j\leq n-1}{\min} V(w_j)>V(w_n),\, V(w_n)\in I(x),\,  \sum_{j=1}^n \underset{u\in \B(w_j)}{\sum}\ee^{-V(u)-x}\leq L  \}}  \{ (1)+ (2) \}  \right)
\end{align*}
with
\begin{align*}
	&(1) = \sum_{k\geq 0}\sum_{v\geq w_n,\, |v|= k+n} 1_{\{  \min_{j\leq k-1} V(v_j)>V(v),\, V(v)\in I(x)   \}},
	\\
&	(2)= \sum_{j=0}^n \sum_{u\in \B(w_j)}\sum_{k\geq 0}\sum_{v\geq u,\, |u|=k+j} 1_{\{  \min_{j\leq k-1} V(v_j)>V(v),\, V(v)\in I(x)   \}}       .
\end{align*}
The term $(1)$ gathers the terms with $v>w_n$ whereas $(2)$ corresponds to those for which $|v\wedge w_n|<|w_n|$. Let $(\tilde{S}_n)_{n\geq 0}$ a independent copy of $(S_n)_{n\geq 0}$.  By the branching property one has
\begin{align*}
(1) &\leq \sum_{n\geq 0} \E\Big(   \1_{\{ \underset{j\leq n-1}{\min} S_j>S_n,\, S_n\in I(x)  \}} \sum_{k\geq 0} \E\big( e^{\tilde{S}_k+a} 1_{\{ \min_{j\leq k-1} \tilde{S}_j> \tilde{S}_k,\, \tilde{S}_k +a\in I(x)  \}} \big)_{a=S_n} \Big)
 \\
 &\leq c e^{-x}\sum_{n\geq 0} \P\left(  \min_{j\leq n-1} S_j>S_n,\, S_n\in I(x)    \right) \leq c'\ee^{-x},
\end{align*} 
where we used twice that
\begin{align}
\label{pascher}
\sup_{x\in \r}\sup_{a\in \r}\sum_{k\geq 0}\P\left(\min_{j\leq n-1} S_j>S_n,\, S_n +a\in I(x)   \right)\leq \sup_{x\in \r}\sup_{a\in \r} [R^{-}(1+x-a)-R^{-}(x-a)] \overset{(\ref{blackwi})}{<} +\infty
\end{align}
To treat the second term, we take the conditional expectation with respect to the sigma-field generated by $(V(w_j),\, (V(u),\, u\in \B(w_j)),\, j\in [|1,n|])$. By the branching property and by recalling that for any $u\in \B(w_j)$,   
$$\E\left(  \sum_{k\geq 0}\sum_{|u|=k} \1_{\{  \min_{j\leq k-1} V(v_j)>V(v),\, V(v) +V(u)\in I(x)   \}}  \right) \leq \ee^{-(x+V(u))}, $$
we get that $ \E\left( N_L(x)^2\right)$  is smaller than
\begin{align*}
	&   \E(N_{L}(x))+  \sum_{n\geq 0} \E_{\Q}\Big( \ee^{V(w_n)}     \1_{\{ \underset{j\leq n-1}{\min} V(w_j)>V(w_n),\, V(w_n)\in I(x),\,   \sum_{j=1}^n \underset{u\in \B(w_j)}{\sum}\ee^{-V(u)-x}\leq L \}}\sum_{j=1}^n\underset{u\in \B(w_j)}{\sum}\ee^{-V	(u)-x}     \Big)  
\\
&\leq c\ee^{-x}+ L \sum_{n\geq 0} \E\left( \ee^{S_n}     \1_{\{ \min_{j\leq n-1} S_j>S_n,\, S_n\in I(x)  \}}  \right)   \leq c' \ee^{-x},
	\end{align*}
where in the last inequality we used (\ref{pascher}). Finally, going back to Paley-Zygmund inequality we have showed that
$$
\P\left(  \mathtt{M}\leq -x \right)\geq \frac{\E( N_L(x) )^2}{\E(N_L(x)^2)}\geq \frac{c'\ee^{-2x}}{c  \ee^{-x}}= c'' \ee^{-x}.$$
It concludes the the proof of the Lemma \ref{tight}. \hfill$\Box$

\subsection{Proof of Theorem \ref{tailMin}}
\label{proffttailmin}
The proof of Theorem \ref{tailMin} requires to study the genealogy of the vertex reaching the global minimum of the branching random walk. This study relies heavily on the {\it spine decomposition} (Proposition \ref{lyons}).  Mention that in the particular case where the displacements have no atom, the global minimum is reached in one unique vertex, what would simplify the computations. 
\\
For any $u\in \mathbb{T}$, $v\in \mathbb{T}^{(u)}$ and $j,k\geq 0,\, a\in \r$, we introduce
\begin{align} 
	\label{Mu0}V^{(u)}(v):=V(v)-V(u),\qquad \M^{(u)}:= \inf\{ V^{(u)}(z),\, z\in \mathbb{T}^{(u)}   \}, \qquad \M^{(u)}_k:= \inf\{ V^{(u)}(z),\, z\in \mathbb{T}^{(u)} \}.
\end{align}
Each depends only on the branching random walk rooted at $u$. Moreover $\M^{(u)}$ and $\M_k^{(u)}$ are distributed respectively as $\M^{(\emptyset)}$ and $\M_k^{(\emptyset)}$ (in the following we shall drop the superscript ${(\emptyset)}$).
\\

\noindent{\it Proof of Theorem \ref{tailMin}.} Recall that $\u$ is chosen at random among the youngest vertices reaching the minimum. Then, by Proposition \ref{lyons}, we get that
\begin{align}
 \nonumber \E\left( \varphi\Big( \mathfrak{D}^{\u,\geq t} \Big)  \1_{\{  \M\in I(x)\}} \right)   & = \sum_{k\geq 0} \E\left(\frac{1}{ N_k(\M)}  \sum_{|z|=k} \1_{\{ V(z)=\M<\M_{k-1},\, V(z)\in I(x)  \}}   \varphi\Big(  \mathfrak{D}^{z,\geq t} \Big) \right)
\\
\label{display1} &= \sum_{k\geq 0} \E_{\Q_k\otimes\P}\left( \frac{\ee^{V(w_k)}}{ N_k(\M)} \1_{\{ V(w_k)=\M<\M_{k-1},\, V(w_k) \in I(x) \}}  \varphi\Big( \mathfrak{D}^{w_k,\geq t} \Big) \right),
\end{align}
with $ \Q_k\otimes\P $ the probability defined in Section \ref{prelimaries} and for any $z\in \T$ with $|z|=k$,
\begin{align*}
\mathfrak{D}^{z,\geq t}:=\sum_{j=k-t}^{k} e^{V(z_k)-V(z_{j-1})} \sum_{v \in \B(z_j)} e^{- [V(v)-V(z_{j-1})] }   D_\infty^{(v)} + D_\infty^{(z_k)}.
\end{align*} 
The event $\{ V(w_k)=\M<\M_{k-1}\}$ can be re-written as
\begin{align*}
 \{  V(w_k)<\min_{j\in [|0,k-1|]}V(w_j)  \}\bigcap \Big(\underset{j\in [|1,k|],u\in \B(w_j)}{\bigcap}\{  V(u)+ \M_{k-j-1}^{(u)}> V(w_k),\,    V(u)+\M^{(u)}\geq V(w_k)  \} \Big)\bigcap \{ \M^{(w_k)}\geq 0 \}
\end{align*}
%Similarly on can decompose $\eta^{\mathtt{M}}(A)$ according to the {\it spine} and write
%\begin{align*}
%\eta^{\mathtt{M}}(A)&= \sum_{v\in \mathbb{T}^{(w_k)}} \1_{\{    \underset{j\in [|0,|v|-1|]}{\min } {V}(v_j)>V(v),\, V(v)-\mathtt{M}\in A   \}}+ \sum_{j=1}^k \sum_{u\in \mathbb B(w_j)} \sum_{v\in \mathbb{T}^{(u)}} \1_{\{    \underset{j\in [|0,|v|-1|]}{\min } {V}(v_j)> V(v),\, V(v)-\mathtt{M}\in A  \}}
%\\
%&\qquad\quad\qquad\qquad\qquad\qquad\qquad\qquad\qquad  \qquad\qquad\qquad   + \sum_{j=1}^k \1_{\{ \min_{l\in [|0,j-1|]}V(w_l)>V(w_j),\, V(w_j)-V(w_k)\in A \}}
%\\
%&= \eta^{(w_k),\mathtt{M}}(A,0,0)+ \sum_{j=0}^k \sum_{u\in \mathbb B(w_j)} \eta^{(u),\mathtt{M}}(A,V(u)-\min_{j\in[0,|u|-1]}V(u_j),V(w_k)-V(u)) +\eta^{\M}(w_k) .
%\end{align*}
On ${\{ V(w_k)=\M<\M_{k-1}  \}}$ one can decompose
\begin{align*}
N_k(\M)=1+\sum_{j=1}^k \sum_{u\in \B(w_j)} \sum_{v\geq u,\, |v|=k}\1_{\{ V(v)=V(w_k) \}}.
\end{align*}
Using these decompositions with (\ref{display1}), by the branching property we get
\begin{align}
\label{display2} \E\left( \varphi\Big(\mathfrak{D}^{\u,\geq t}   \Big)  \1_{\{  \M\in I(x)\}} \right) = \sum_{k\geq 0} \E_{\Q_k\otimes \P}\left( \frac{1}{N_k(\M)} { \ee^{V(w_k)}}  \1_{\{  V(w_k)<\underset{j\in [|0,k-1|]}{\min}V(w_j),\, V(w_k)\in I(x) \}}    \right.
\\
  \prod_{j=1}^{k} \prod_{u\in \mathbb{B}(w_j)}   1_{\{ \mathtt{M}^{(u)}\geq  V(w_k)-V(u),\, \mathtt{M}_{k-j-1}^{(u)}>  V(w_k)-V(u)  \}}  \varphi\Big(\mathfrak{D}^{w_k,\geq t}  \Big)      \Big).\nonumber 
\end{align}
When $x$ goes to $+\infty$ many terms of this expression can be simplified. Indeed the following two lemmas will state that the first terms of the infinite sum above are negligible and that all the particles of the branching random walk whose the position is close to $\mathtt{M}$ are also genealogically close to $w_k$.
\begin{Lemma}
	\label{Cuty} For any $b\geq 0$, 
	\begin{equation}
		\label{lefthand}
	\lim_{x\to\infty}  \ee^x	\sum_{k\leq b}\E_\Q\left( \ee^{V(w_k)}\1_{\{   V(w_k)<\min_{j\in [|0,k-1|]}V(w_j),\, V(w_k)\in I(x) \}} \right) =0.
	\end{equation}
\end{Lemma}
Under the probability $\Q_k\otimes \P$,  let us define
 \begin{align}
 \mathcal{E}_k(b_1)&:= \{\forall j\leq k-b_1,\, \forall u\in \mathbb B({w_j}),\, V(u)+ \mathtt{M^{(u)}}\geq  V(w_k)+ 1 \}, \quad\, \forall b_1\geq 0
 %\\
 %F_k( b_1, b_2)&:= \{ \max_{j\in [k-b_1,k]}V(w_j) <\min_{j\in [ 0,k-b_2]}V(w_j) \},\qquad\quad \qquad \forall  b_2\geq b_1\geq 0.
 \end{align}
\begin{Lemma}
	\label{cuty2}For any $D>0$,  
	\begin{equation*}
\lim_{b_1\to\infty} \lim_{b_2\to\infty}\lim_{x\to \infty} \ee^x	\sum_{k\geq b_2}\E_{\Q_k\otimes\P}\left( \ee^{V(w_k)}\1_{\{   V(w_k)<\min_{j\in [|0,k-1|]}V(w_j),\, V(w_k)\in I(x) \}}  \1_{ \mathcal{E}_k(b_1)^c  }   \right)=0
	\end{equation*}
%	and
%	\begin{equation}
%	\label{BIS3.3}
%	\lim_{b_1\to\infty} \lim_{b_2\to\infty}\lim_{x\to \infty} \ee^x	\sum_{k\geq b_2}\E_{\Q_k\otimes\P}\left( \ee^{V(w_k)}\1_{\{    V(w_k)<\min_{j\in [|0,k-1|]}V(w_j),\, V(w_k)\in I(x) \}}   \1_{  F_k(b_1,b_2)^c  }    \right)=0
%	\end{equation}
\end{Lemma}
The proofs are postponed in the next subsection. Applying Lemma \ref{Cuty} and \ref{cuty2}, we can affirm that: {\it For any $\epsilon >0$, there exists $x_0>0$ such that for any $x\geq x_0$ there exits $B_1>0$ such that for any $b_1\geq B_1$, there exists $B_2>b_1$ such that for any $b_2\geq B_2$}
\begin{align}
\label{bious}
&\Big|\E\Big(\varphi\Big( \mathfrak{D}^{\u,\geq t}  \Big)    1_{\{ \M\in I(x) \}}    \Big)-    \sum_{k\geq b_2}  \E_{(\ref{bious})}^{(k)}(b_1,b_2) \Big|\leq \epsilon \ee^{-x}
\end{align}
with  for any $k\geq b_2  \geq b_1 $, 
\begin{align*} 
& \E_{(\ref{bious})}^{(k)}(b_1,b_2):= \E_{\Q}\Bigg( \frac{  { \ee^{V(w_k)}}  \1_{\{  V(w_k)<\underset{j\in [|0,k-1|]}{\min}V(w_j),\, V(w_k)\in I(x) \}}   }  {  1+\sum_{j=k-b_1+1}^k \sum_{u\in \B(w_j)} \sum_{v\geq u,\, |v|=k}\1_{\{ V(v)=V(w_k) \}}}     
\\
& \qquad\qquad \qquad\qquad  \prod_{j=k-b_1+1}^{k} \prod_{u\in \mathbb{B}(w_j)}   1_{\{ \mathtt{M}^{(u)}\geq  V(w_k)-V(u),\, \mathtt{M}_{k-j-1}^{(u)}>  V(w_k)-V(u)   \}}   \varphi\Big(  \mathfrak{D}^{w_k,\geq t} \Big)      \Bigg).
%\\
%&=\E_{\Q_k\otimes \P}\left( \ee^{V(w_k)}\1_{\{  V(w_k)=\underline{V}(w_k),\, V(w_k)\in I(x) \}} \times\right.
%\\
%& \qquad\qquad \prod_{j=k-b_1}^{k} \prod_{u\in \mathbb{B}(w_j)}  \E_\P\left(  1_{\{ \mathtt{M} \geq x_1  \}} \ee^{-\theta \eta^{ \mathtt{M}}(A,x_2, x_3) }   \right)_{\Big| x_1= V(w_k)-V(u), x_2=V(u)-\underset{l\in [k-b_2,j-1]}{\min}{V}(u_l), x_3=V(w_k)-V(u) }  \Big)
%\\
%&=\E_{\Q_k\otimes \P}\left( \ee^{V(w_k)}\1_{\{  V(w_k)=\underline{V}(w_k),\, V(w_k)\in I(x) \}} \times\right.
%\\
%& \qquad\qquad \prod_{j=k-b_1}^{k} \prod_{u\in \mathbb{B}(w_j)}  \tilde{\Psi}_A^{(\theta)}(  V(w_k)-V(u),  V(u)-\underset{l\in [k-b_2,j-1]}{\min}{V}(u_l),  V(w_k)-V(u) )  \Big).
\end{align*}
%with 
%$$ \eta_{b_1,b_2}^{\M,k}:= \sum_{j=k-b_1}^k \1_{\{ \min_{l\in [|k-b_2,j-1|]}V(w_l)>V(w_j),\, V(w_j)-V(w_k)\in A \}} $$
For any $j\in [|1,k|]$, recall that $\sum_{u\in \B (w_j)}\delta_{\zeta^{(j)}_u }=  \sum_{u\in \B (w_j)}\delta_{V(u)-V(w_{j-1}) }$ is the point process formed by the brothers of $w_j$, and for any $j\in \N^*$ and $u\in \B(w_{j})$, let $
\Xi^{(j)}_u:= V(w_{j})-\zeta_u^{(j)}$. By using the time reversal identity (\ref{reversal}) we obtain the following changes:
\begin{align*}
V(w_k) \leftrightarrow V(w_k),\qquad V(w_k) <  \min_{j\in [|0,k-1|]}V(w_j) \leftrightarrow  \max_{j\in [|1,k|]}V(w_j) <0,\qquad u\in \B(w_j) \leftrightarrow u\in \B(w_{k-j+1}).
\end{align*}
When $u\in \B(w_j)$:
\begin{align*}
  V(w_k)-V(u) = V(w_k) -V(w_{j-1}) - \zeta_u^{(j)} &\leftrightarrow  V(w_{k-j+1}) -\zeta_u^{(k-j+1)}= \Xi_u^{(k-j+1) },\qquad \text{with} \quad u\in B(w_{k-j+1}),
\\
 V(u) - \min_{l\in [|k-b_2,j-1|]}V(u_l)  &\leftrightarrow     \Xi_u^{(k-j+1)} +  \max_{l\in [|k-j+1, b_2  |]} V(w_l),\qquad\qquad   \text{with} \quad u\in B(w_{k-j+1}),
\end{align*}
and
\begin{align*}
\mathfrak{D}^{\u,\geq t} &\leftrightarrow \sum_{j=k-t}^{k} e^{V(w_{k-j+1})} \sum_{v \in \B(w_{k-j+1})} e^{-  \zeta_{v}^{(k-j+1)} }   D_\infty^{(v)} + D_\infty^{(w_0)} 
\\
&=  \sum_{j=1}^{t+1} e^{V(w_{j})} \sum_{v \in \B(w_{j})} e^{-  \zeta_{v}^{(j)} }   D_\infty^{(v)} + D_\infty^{(w_0)}  :=\bar{\mathfrak{D}}^{w_{t+1},\leq t+1},
\end{align*}
where we recall that $(D_\infty^{(v)})_{v\in \B(w_j),j\in [|1,t+1|]}$, $D_\infty^{(w_0)}$ are the limit of the {\it Derivative martingales} of the branching random walks rooted respectively at $v$ and $w_0$.  Finally one can write
\begin{align*}
&\E_{(\ref{bious})}^{(k)}(b_1,b_2)=\E_\Q\left(   \frac{ \ee^{V(w_k)}\1_{\{ \max_{j\in [|1,k|]}V(w_j) <0,\,   V(w_k)  \in I(x)     \}}  }{1+  \sum_{j=k-b_1+1}^k \sum_{u\in \B(w_{k-j+1})} \sum_{ |v|=k-j+1}\1_{\{   V^{(u)}(v)=  \Xi_u^{(k-j+1)} \}}  }  \right.
\\
& \qquad\qquad\qquad\qquad\qquad  \left.\prod_{j=k-b_1+1}^k \prod_{u\in  \B(w_{k-j+1})} \1_{\{\M^{(u)}\geq \Xi_u^{(k-j+1)} ,\, \M_{k-j}^{(u)}>  \Xi_u^{(k-j+1)}  \}}  \varphi\Big(  \bar{\mathfrak{D}}^{w_{t+1},\leq t+1} \Big)  \right)
\end{align*}
Now by operating the change of index $j\leftrightarrow k-j+1$ in the product and in the denominator, it becomes
\begin{align*}
\E_{(\ref{bious})}^{(k)}(b_1,b_2)
&=\E_\Q\left(   \frac{ \ee^{V(w_k)}\1_{\{ \max_{j\in [|1,k|]}V(w_j)<0,\, V(w_k) \in I(x)     \}}     }{1+\sum_{j=1}^{b_1} \underset{u\in \B(w_{j})}{\sum} \underset{ |v|=j}{\sum}\1_{\{ V^{(u)}(v)=  \Xi_u^{(j)} \}}  }     \prod_{j=1}^{b_1} \prod_{u\in  \B(w_{j})} \1_{\{\M^{(u)}\geq \Xi_u^{(j)} ,\, \M_{j-1}^{(u)} >  \Xi_u^{(j)}  \}} \varphi\Big(  \bar{\mathfrak{D}}^{w_{t+1},\leq t+1}  \Big)  \right).
\end{align*}
By applying the branching property at the vertex $w_{b_2}$, $ \ee^x\sum_{k\geq b_2}\E^{(k)}_{(\ref{bious})}(b_1,b_2)$ is equal to
\begin{align}
\label{imingbac0} &   \E_\Q\left(  \frac{  \1_{\{ \underset{j\in [|1,b_2|]}{\max}   V(w_j)<0\}}   \E_{V(w_{b_2})}\left( \sum_{k\geq 0}\ee^{S_k+x}\1_{\{   \underset{j\in [|1,k|]}{\max}S_j<0,\, S_k\in I(x)   \}} \right)}{1+\sum_{j=1}^{b_1} \underset{u\in \B(w_{j})}{\sum} \underset{ |v|=j}{\sum}\1_{\{ V^{(u)}(v)=  \Xi_u^{(j)} \}}  } \varphi\Big(  \bar{\mathfrak{D}}^{w_{t+1},\leq t+1}  \Big) \prod_{j=1}^{b_1} \prod_{u\in \B(w_{j})}... \right),
%
%   \right.
%\\
%\nonumber& \left. \times \prod_{j=0}^{b_1} \prod_{u\in \B(w_{j+1} )} \1_{\{\M^{(u)}\geq  V(w_{j+1})-\zeta_u^{(j+1)} ,\, \M_{j-1}^{(u)} >  V(w_{j+1})-\zeta_u^{(j+1)}  \}}\ee^{-\theta \eta^{(u),\M} (A, \zeta_u^{(j+1)}- V(w_{j+1}) +\underset{l\in [|j+1,b_2|]}{\max}V(w_l), V(w_{j+1})-\zeta_u^{(j+1)} )  }     \right) .
\end{align}
Note that only the $\E_{V(w_{b_2})}(...)$ term depends on the variable $x$. Furthermore by standard computations, for any $x,a\geq 0$, 
\begin{align}
\nonumber  \E_{-a}\left( \sum_{k\geq 0} \ee^{S_k +x}\1_{\{  \underset{j\in [|1,k|]}{\max}S_j<0,\,  S_k\in I(x)  \}}   \right)  &= \ee^{x}\sum_{k\geq 0} \E_{-a}\left( \int_{-\infty}^{+\infty} \ee^{u} \1_{\{u\leq  S_k \}} du\,  \1_{\{   \max_{j\in [|1,k|]}S_j< 0,\, S_k\in I(x)  \}}   \right)
 \\
\nonumber  &= \ee^x \int_{-\infty}^{-x} \ee^{u} \sum_{k\geq 0} \P_{-a} \left(  \max_{j\in [|1,k|]}S_j<0,\, -x> S_k \geq \max(-x-1,u)\right) du
 \\
\label{dComoi}\nonumber  &= \ee^x \int_{-\infty}^{-x} \ee^u [  \tilde{R}(\min(x+1,-u),a)-\tilde{R}(x,a)   ]du
\\
&= \ee^{-1} [\tilde{R}(x+1,a)-\tilde{R}(x,a)] + \int_{0}^1 \ee^{-u}[\tilde{R}(x+u,a)-\tilde{R}(x,a)]du ,
\end{align}
where $\tilde{R}$ is the function defined in (\ref{tildiR}). Using (\ref{chap3lldda}), Lemma  \ref{renouvDoux} and (\ref{laRemark}), it follows that
\begin{align*}
\lim_{x\to\infty} \E_{-a}\left( \sum_{k\geq 0} \ee^{S_k +x}\1_{\{  \underset{1\leq j\leq k}{\max}S_j< 0,\,  S_k\in I(x)  \}}   \right)&= C^- \theta_0( R^+(a)-K_a) \left\{ \ee^{-1} + \int_0^1 u\ee^{-u}du    \right\}
\\
&= (1-\ee^{-1})C^-  \theta_0 ( R^+(a)-K_a )  .
%\\
%& =\theta_0\{R^+(a)-K_a\} \int_{-\infty}^{-x}\ee^{u+x} \{ R^-( \min(x+1,-u))-R^-( x)  \}du??
%\\
%&= \theta_0\{R^+(a)-K_a\} \int_{0}^{+\infty}\ee^{-u} \{ R^-( \min(x+1,x+u))-R^-( x)  \}du
%\\
%&= \theta_0\{R^+(a)-K_a\} C^-  \left( \int_{0}^1 u \ee^{-u} du + \int_1^{+\infty} \ee^{-u} du   \right)
%\\
%&= c^+\{R^+(a)-\P(\exists i\in \N, S_i=a)\} C^-(1-\ee^{-1}).
\end{align*}
%Finally by letting $x$ going to infinity we have
%\begin{equation}
%\lim_{x\to\infty}  \E_{-a}\left( \sum_{k\geq 0} \ee^{S_k +x}\1_{\{ \underset{j\in [|1,k|]}{\max} S_j< 0,\,  S_k\in I(x)  \}}   \right)= c^+R^+(a)C^-={\red K_0\{R^+(a)-\P(\exists i\in \N, S_i=a)\}  }.
%\end{equation}
Plugging this equality in (\ref{imingbac0}) we have
  \begin{align*}
&\lim_{x\to \infty} \ee^x\sum_{k\geq b_2}\E^{(k)}_{(\ref{bious})}(b_1,b_2)=
\\
& (1-\ee^{-1})C^-\theta_0   \E_\Q\left(  \frac{   ( R^+(-V(w_{b_2}))-K_{-V(w_{b_2})} )  \1_{\{ \underset{j\in [|1,b_2|]}{\max}V(w_j)<0\}} }{  1+\sum_{j=1}^{b_1} \underset{u\in \B(w_{j})}{\sum} \underset{ |v|=j}{\sum}\1_{\{ V^{(u)}(v)=  \Xi_u^{(j)} \}}  } \varphi\Big(  \bar{\mathfrak{D}}^{w_{t+1},\leq t+1}  \Big)  \prod_{j=1}^{b_1} \prod_{u\in \B(w_{j})} \1_{\{\M^{(u)}\geq \Xi_u^{(j)} ,\, \M_{j-1}^{(u)} >  \Xi_u^{(j)}  \}}\right).
%  \right.
%\\
%& \left. \times \prod_{j=0}^{b_1} \prod_{u\in \B(w_{j+1})} \1_{\{\M^{(u)}\geq V(w_{j+1})-\zeta_u^{(j+1)} ,\, \M_{j-1}^{(u)} > V(w_{j+1})-\zeta_u^{(j+1)}  \}}\ee^{-\theta \eta^{(u),\M} (A, \zeta_u^{(j+1)}- V(w_{j+1}) +\underset{l\in [|j+1,b_2|]}{\max}V(w_l), V(w_{j+1})-\zeta_u^{(j+1)} )  }     \right).
\end{align*}
When $b_2$ goes to infinity the term induced by $K_{-V(w_{b_2})}\leq 1$ converges to $0$. Then by the monotonicity of $\sum_{k\geq b_2}\E^{(k)}_{(\ref{bious})}(b_1,b_2)$; in $b_1$ and $b_2$, we deduce that the following limit exists
\begin{align*}
\nonumber &\lim_{x\to \infty} \ee^{x} \E\left(  \varphi\Big(\mathfrak{D}^{\u,\geq t}   \Big)\1_{\{  \M\in I(x)\}} \right)= (1-\ee^{-1})   \times \mathcal{E}_t(\varphi),
\end{align*}
with  $\mathcal{E}_t(\varphi)$ defined by
\begin{align}
\label{defunctio}
C^-\theta_0\lim_{b_1\to\infty} \lim_{b_2\to\infty}  \E_{  \Q}\left(  \frac{  \varphi\Big( \bar{\mathfrak{D}}^{w_{t+1},\leq t+1}  \Big)  R^+(-V(w_{b_2}))   \1_{\{ \underset{j\in [|1,b_2|]}{\max}V(w_j)<0\}} }{  1+\sum_{j=1}^{b_1} \underset{u\in \B(w_{j})}{\sum} \underset{ |v|=j}{\sum}\1_{\{ V^{(u)}(v)=  \Xi_u^{(j)} \}}  }   \prod_{j=1}^{b_1} \prod_{u\in \B(w_{j})} \1_{\{\M^{(u)}\geq  \Xi_u^{(j+1)} ,\, \M_{j-1}^{(u)} >\Xi_u^{(j+1)}   \}}     \right).
\end{align}
By using this convergence and replacing $I(x)$ by  $I(x+1)$, $I(x+2)$, $I(x+3)$... and summing everything we obtain Theorem \ref{tailMin}.  \hfill$\Box$

%\begin{align}
%\lim_{x\to\infty}  	\sum_{k\geq b_2}   \E_{(\ref{bious})}^{(k)}(b_1,b_2)&= \E\left(  \1_{\{ \overline{S}_{b_2}\leq 0 \}} \times \prod_{j=0}^{b_1} \psi_{A}(  X_j, S_{j+1} ,    \max_{l\in [j+1,b_2]}S_l  -S_{j+1}        ) C^-R^+(S_{b_2}) \right)\nonumber
%\\
%&= C^-\E_0\left(   \prod_{j=0}^{b_1} \psi_{A}(  X_j, S_{j+1} ,    \max_{l\in [j+1,b_2]}S_l  -S_{j+1}        )   \right),
%\end{align}
%where under $\E_0$ the process $(S_n)_{n\geq 0}$ is distributed as a random walk conditioned to stay negative. This last expression is monotonous in $b_1$ and $b_2$, then letting $b_1$ then $b_2$ going to infinity we obtained that 
%{\green \begin{align}
%\lim_{x\to \infty} \ee^{x} \E\left( \exp( -\theta \eta^{\M}(A)) \1_{\{  \M\leq -x\}} \right)=  C^-\E_0\left(   \prod_{j=0}^{\infty} \psi_{A}^{(\theta)}(  X_j, S_{j+1} ,    \max_{l\in [j+1,+\infty]}S_l  -S_{j+1}        )   \right).
%\end{align}}
%
% 
% 

\subsection{Proof of Lemma \ref{Cuty} and Lemma \ref{cuty2}} 

\noindent{\it Proof of Lemma \ref{Cuty}.} Fix $b>0$. The left-hand term of (\ref{lefthand}) is equal to 
\begin{align*}
&\sum_{k\leq b} \ee^x\E\left( \ee^{S_k}1_{\{ S_k<   \underset{j\in [|0,k-1|]}{\min } {S}_j,\, S_k\in [-x-1,-x)  \}}\right) \leq \sum_{k\leq b} \P\left(  S_k\in [-x-1,-x]  \right) \to 0,\qquad \text{when }x\to \infty,
\end{align*}
which concludes the proof of Lemma \ref{Cuty}.
\hfill$\Box$ 

\noindent{\it Proof of Lemma \ref{cuty2}.} Recall that we need to prove  that
	\begin{equation}
		\label{bisFisrt3}
 	\lim_{b_2\to\infty} \lim_{b_1\to\infty}\lim_{x\to \infty}  	\sum_{k\geq b_2}\E_{\Q_k\otimes\P}\left( \ee^{x+ V(w_k)}\1_{\{  V(w_k)< \underset{j\in [|0,k-1|]}{\min } {V}(w_j),\, V(w_k)\in I(x) \}} ; \mathcal{E}_k(b_1)^c\right)=0.
	\end{equation}
Let us denote $\E^{(k)}_{(\ref{bisFisrt3})}$ the expectation in (\ref{bisFisrt3}). Notice that
\begin{align}
\E^{(k)}_{(\ref{bisFisrt3})}\leq \sum_{j=1}^{k-b_1}  {\Q_k\otimes\P}\left(    V(w_k)< \underset{j\in [|0,k-1|]}{\min } {V}(w_j),\, V(w_k)\in I(x),\,  \exists u\in \B(w_j),\, V(u) +\M^{(u)}\leq V(w_k) +1   \right) .
\end{align}
Moreover for any $k\geq b_2\geq b_1$, $j\leq k-b_1$, by the branching property and (\ref{eqticht}) we have
\begin{align*}
&\Q_k\otimes\P\left(  \exists u\in \B(w_j),\, V(u) +\M^{(u)}\leq V(w_k) +1 \big| \sigma \left( (V(w_l), V(u),\, u\in \B(w_l))_{l\in [|1,k|]}  \right)\right) 
\\
  &\qquad \qquad\qquad \leq  \min( 1,\sum_{u\in \mathbb B(w_j)}\ee^{V(w_k)-V(u) +1} ) =  \min(1,  \ee^{V(w_k)-V(w_{j-1}) +1}\Delta_j),
\end{align*}
%
%
%\begin{align*}
%&\sum_{k\geq b_2} \ee^x\E_{\Q_k\otimes\P}\left( \ee^{V(w_k)}\1_{\{  V(w_k)<   \underset{j\in [|0,k-1|]}{\min } {V}(w_j),\, V(w_k)\in I(x) \}} \sum_{j=1}^{k-b_1} \1_{\{   \sum_{u\in \mathbb B(w_j)}\1_{\{ V(u)+ \M^{(u)}\leq -x +D  \}}  >0\}} \right)
%\\
%&\leq \sum_{k\geq b_2} \ee^x\E_{\Q_k\otimes\P}\left( \ee^{V(w_k)}\1_{\{ V(w_k)<   \underset{j\in [|0,k-1|]}{\min } {V}(w_j)  ,\, V(w_k)\in I(x) \}} \sum_{j=1}^{k-b_1}\min( 1,\sum_{u\in \mathbb B(w_j)}\ee^{-(x+V(u)+ D)} )\right)
%\end{align*}
%Then we have
%\begin{align*}
%&\leq \sum_{k\geq b_2} \sum_{j=0}^{(k-b_1)_+} \E_{\Q}\left( \1_{\{  V(w_k)=   \underset{j\in [|0,k|]}{\min } {V}(w_j),\, V(w_k)\in I(x) \}} \min[ 1, \ee^D \ee^{V(w_k)-V(w_j)} \Delta_j]\right)
%\end{align*}
with $\Delta_j:= \sum_{u\in \B(w_j)}e^{V(w_{j-1})-V(u)}$. Recall that $(V(w_j)-V(w_{j-1}),\Delta_j)_{j\in [|1,k|]}$ are i.i.d., then by operating a time reversal we have
\begin{align*}
 \sum_{k\geq b_2} \E^{(k)}_{(\ref{bisFisrt3})}   &\leq \sum_{k\geq b_2} \sum_{j=1}^{k-b_1}   \E_\Q \left(  \min [1, \ee^1\Delta_{k-j+1} \ee^{ V(w_{k-j+1})} ]\1_{\{ \underset{l\in [|1,k|]}{\max}V(w_l)< 0,\, V(w_k)\in I(x)  \}} \right)
\\
&\leq c \sum_{k\geq b_2} \sum_{j=b_1+1}^{k}   \E_\Q \left(  \min [1, \Delta_{j} \ee^{ V(w_{j})} ]  \1_{\{ \underset{l\in [|1,k|]}{\max}V(w_l)< 0,\, V(w_k)\in I(x)  \}}  \right).
\end{align*}
By the branching property at time $j$, for any $x\geq 1$ we get
%
%\begin{align}
%\sum_{j\in \N} \sum_{k\in \N}  1_{\{  k\geq b_2,\, k\geq j,\, j \geq b_1   \}}= \sum_{j\in \N} \sum_{k\geq \max(j,b_2)}^{+\infty}
%\end{align}
%
\begin{align*}
 \sum_{k\geq b_2} \E^{(k)}_{(\ref{bisFisrt3})} &\leq c\sum_{k\geq b_2} \sum_{j=b_1+1}^{k}   \E_\Q \left(  \min [1,  \Delta_j \ee^{ V(w_{j})} ] 1_{\{    \underset{i\leq j}{\max}V(w_i) <0  \}} \P_{V(w_j)}  \left( \underset{l\in [|1,k-j|]}{\max} S_l < 0,\, S_{k-j}\in I(x)     \right)  \right)
\\ 
&\leq c' \sum_{j=b_1+1}^\infty \E_\Q\left(  \min [1,  \Delta_j \ee^{ V(w_{j})} ]  1_{\{    \underset{i\leq j}{\max}V(w_i) <0  \}}   \sum_{k= 0}^{+\infty} \P_{V(w_j)}  \left(  \underset{i\in [|1,k|] }{\max}\, S_i<0,\, S_k\in I(x)     \right) \right)
\\
&= c' \sum_{j=b_1+1}^\infty \E_\Q\left(  \min [1,  \Delta_j \ee^{ V(w_{j})} ]  1_{\{    \underset{i\leq j}{\max}V(w_i) <0  \}} \left[ \tilde{R}(x+1,-V(w_j))  -\tilde{R}(x,-V(w_j))   \right] \right)
\\
&\leq   c''\sum_{j=b_1+1}^\infty \E_\Q\left(  \min [1,  \Delta_j \ee^{ V(w_{j})} ]  \1_{\{    \underset{i\leq j}{\max}V(w_i) <0  \}} (1-V(w_j) ) \right),
\end{align*}
where in the last inequality we used (\ref{consequencese}). Notice that the last expression does not depend in $x$ any more. Moreover for any $L>0$, 
\begin{align*}
\min[1,\Delta_j\ee^{V(w_j)}] &\leq \ee^{L+ \frac{1}{2}V(w_j)} + \1_{\{ \ln \Delta_j \geq   L-\frac{1}{2}V(w_j) \}}
\\
&\leq \ee^{L+ \frac{1}{2}V(w_j)} + \1_{\{-V(w_{j-1}) \leq V(w_j)-V(w_{j-1})+ 2\ln \Delta_j-2L \}}
\end{align*}
Let $(\Delta,\zeta)$ a couple of random variables distributed as $(\Delta_1, V(w_1))$ and independent of everything else. By using the inequality just above, for any $b_1,L>0$ one has
\begin{align*}
 \sum_{k\geq b_2} \E^{(k)}_{(\ref{bisFisrt3})} &\leq c\sum_{j=b_1}^\infty \E_{\Q}\left( (1-S_j )  \ee^{L+\frac{1}{2}S_j}\1_{\{    \underset{l\in [|0,j|]}{\max } S_l\leq 0 \}}\right)  +c\sum_{j=b_1}^\infty \E_{\Q}\left((1+ \ln \Delta)  \1_{\{  -S_{j-1}\leq \zeta+  2 \ln \Delta -2L \}}  \1_{\{  \underset{l\in [|0,j-1|]}{\max } S_l\leq 0\}} \right) 
 \\
 &\leq c'\ee^L\sum_{j=b_1}^\infty \E_{\Q}\left( \ee^{\frac{1}{4}S_j}\1_{\{    \underset{l\in [|0,j|]}{\max } S_l\leq 0 \}}\right)  + c'\E_{\Q}\left((1+ \ln\Delta)(1+\ln \Delta+\zeta)\1_{\{ 2L\leq \ln \Delta +\zeta  \}}  \right) .
\end{align*}
 For any $L>0$, when $b_1$ goes to infinity, the first term converges to $0$ by Lemma \ref{techniqAID}. The second one converges to $0$ when $L$ goes to infinity because of the assumption (\ref{chap31.5}) (notice that $\ln \Delta +\zeta$ are stochastically dominated by the random variable $X$) . It concludes the proof of (\ref{bisFisrt3}).\hfill$\Box$
\\

\section{Proof of Proposition \ref{belowplus}} 
Recall that $\mathfrak{D}^\M$ is defined in (\ref{deffrakD}).  The following Lemma studies the integrability of $\mathfrak{D}^\M$. We introduce
\begin{align*}
D_n^{(a)}:= \sum_{|z|=n} R^-(V(z)+a) e^{-V(z)} 1_{\{ \min_{j\in[0,,n]}V(z_j)\geq -a  \}},\qquad n\in \N^*,\, a>0.
\end{align*}
It is a non-negative martingale with mean $a$. In \cite{BKy04}, Biggins and Kyprianou proved that there exists $c_0>0$ such that for any $a>0$,  on $\{\M \geq -a\}$, 
\begin{align}
\label{BigKYP}
\lim_{n\to\infty} D_n^{(a)}= c_0 D_\infty. 
\end{align}
Recall also that $a>0$,  
\begin{align}
\label{piousit}
\E\left( D_\infty \1_{\{ \M\geq -a\}}\right) \leq  \frac{1}{c_0}\E\left( D_\infty^{(a)}  \right) = \frac{1}{c_0}a.
\end{align}

\noindent{\it Proof of Proposition \ref{belowplus}.}  Recall that $I(x)=[-x-1,x)$. Let us define 
  \begin{align}
  h(u)= \left\{  \begin{array}{ll} u,\qquad  &\text{for } \,\, u<\ee^1,
  \\
  \ee^{1} -1+ \ln^2(u), &\text{for } \,\,u\geq \ee^1.
  \end{array} \right.
  \end{align}
  It is plain to check that $h$ is concave, continuous and increasing. Moreover inequality (\ref{below2plus}) is equivalent to the following inequality  
  \begin{align*}
   	\sup_{x\in \r^+} \E\left(     \mathfrak{D}_x^\M   h(  \mathfrak{D}_x^\M   ) \right)<\infty.
  \end{align*} 
 By Lemma \ref{tight}, and (\ref{BigKYP}), 
 \begin{align*}
 \sup_{x\in \r^+} \E\left(    \mathfrak{D}_x^\M   h(  \mathfrak{D}_x^\M ) \right) \leq c \sup_{x\in \r^+}  e^x \E\left( e^{\M} D_\infty      h(    e^{\M} D_\infty  )  1_{\{  \M\leq -x \} } \right)& \leq c \sup_{x\in \r^+}   \sum_{p\geq 0} e^{-p }   e^{x+p} \E\left(   e^{\M} D_\infty      h(    e^{\M} D_\infty  )  1_{\{ \M\in I(x+p) \}}   \right)
 \\
 &\leq c'\sup_{x\in \r^+}  \sum_{p\geq 0} e^{-p }      \E\left( \lim_{n\to\infty}    D_n^{(x+p+1)}     h(    e^{\M} D_\infty  )  1_{\{ \M \in I(x+p) \}}   \right)
 \\
 &\leq c'\sup_{x\in \r^+}  \sum_{p\geq 0} e^{-p}   \sup_{n\geq 1} \E\left(     D_n^{(x+p+1)}     h(    e^{\M} D_\infty  )  1_{\{ \M\in I(x+p+1) \}}   \right)
 \\
 &\leq  c'' \sup_{x\in \r^+} \sup_{n\geq 1} e^x \E\left(    D_n^{(x+1)}     h(    e^{\M} D_\infty  )   1_{\{ \M\in I(x) \}}   \right).
 \end{align*}
 Finally  it boils down to prove that there exists $c>0$ such that for any $x\geq 0,\, n\geq 1$
 \begin{align}
\nonumber  M_n^h(x)&:= \E\left(   \sum_{|z|=n} R^-(V(z)+x+1) e^{-V(z)} \1_{\{ \min_{j\in[0,,n]}V(z_j)\geq -(1+x+1) \}}  h(    e^{\M} D_\infty  )   \1_{\{ \M  \in I(x)  \}} \right) 
 \\
 \label{plusforloweretoile} &=  \E_\Q\left( R^-(x+1+V(w_n)) 1_{\{ \min_{j\in[0,n]}V(w_j) \geq -x-1  \}} h(e^{\M} D_\infty)   1_{\{ \M \in I(x)  \}}  \right)\leq c,
 \end{align}
where the last equality is justified by Proposition \ref{lyons}. On $\{ \M\in I(x)\}$, one has
  \begin{align*}
  e^{\M} D_\infty &\leq  \sum_{k=1}^{n} \sum_{v \in \B(w_k)} e^{-x -V(v) }   \bar{D}_\infty^{(v,x)} + e^{-x-V(w_n)}\bar{D}_\infty^{(w_n,x)} ,
  \end{align*}
  with $\bar{D}_\infty^{(v,x)}:= \lim_{n\to\infty} \sum_{|z|=n,\, z\geq v} [V(z)-V(v)]e^{-V(z)-V(v)} \1_{\{ \min_{i\in [|v|,n]}V(z_i)\geq -x-V(v)  \}}$. Let us denote by $\sigma(\B_n)$ the sigma field generated by $ (  V(w_j), (V(u),\, u\in \B(w_j)))_{j\in [1,n]}$. By Jensen inequality combined with (\ref{piousit}) we have that
  \begin{align*}
  \E\left( h( e^{\M} D_\infty )\Big| \sigma(\B_n) \right) &\leq    h( c\sum_{k=1}^{n} \sum_{v \in \B(w_k)} (x+1 +V(v))e^{x -V(v) }   + c (x+1+ V(w_n) )e^{x-V(w_n)}   ) 
  \\
  &\leq h( c\sum_{k=1}^{n}  ([x+ V(w_{k-1})]_+ +1)e^{-x-V(w_{k-1})} \tilde{\Delta}_k   + c (x+1+ V(w_n) )e^{x-V(w_n)}   ) 
  \\
  &\leq h(  c' \sum_{k=1}^{n}e^{-\frac{1}{2}(x+V(w_{k-1}))} \tilde{\Delta}_k  +c(x+1+ V(w_n) )e^{x-V(w_n)}  ).
  \end{align*}
  with $\tilde{\Delta}_k:=   \sum_{v \in \B(w_k)} e^{-[V(v)-V(w_{k-1})]} ( [V(v)-V(w_{k-1})]_+ +1)$. Moreover for any $x,y\geq 0$, $h(x+y)\leq h(x)+ h(y)$, it follows that
 \begin{align*}
 M_n^h(x)\leq \sum_{k=1}^n  \E_\Q\left( R^-(x+1+V(w_n)) 1_{\{ \min_{j\in[0,n]}V(w_j) \geq -x-1  \}} h(  c e^{-\frac{1}{2}(x+V(w_{k-1}))} \tilde{\Delta}_k)  \right)
 \\
 +  \E_\Q\left( R^-(x+1+V(w_n)) 1_{\{ \min_{j\in[0,n]}V(w_j) \geq -x-1  \}} h( ce^{- (x+V(w_n))}  )  \right).
\end{align*}
As $h(x)\leq x$ and $R^{-}(x)\leq c(1+x)$ for any $x\geq 0$, the second term is trivially bounded by
\begin{align*}
c\E_\Q\left( R^-(x+1+V(w_n)) 1_{\{ \min_{j\in[0,n]}V(w_j) \geq -x-1  \}}  e^{- (x+V(w_n))}   \right)\leq c'.
\end{align*}
Concerning the sum, after using the Markov property at time $k$ we need to prove that there exists $c>0$ such that for any $x\geq 1$ and $n\in \N$,
\begin{align}
\label{aprouve}\sum_{k=1}^n I_k:= \sum_{k=1}^n  \E_\Q\left( R^-(x+1+V(w_k)) 1_{\{ \min_{j\in[0,k]}V(w_j) \geq -x-1  \}} h( e^{-\frac{1}{2}(x+V(w_{k-1}))} \tilde{\Delta}_k)  \right)\leq c.
\end{align}
By partitioning the expectation on  $\cup_{p\in\N^*} \{ x+V(w_{k-1}) +p   \leq  4\ln \tilde{\Delta_k} \leq x+V(w_{k-1}) +p +1 \}$, one obtains that
\begin{align*}
&\sum_{k=1}^n I_k\leq c\sum_{k=1}^n \E_\Q\left( R^-(x+1+V(w_k)) 1_{\{ \min_{j\in[0,k]}V(w_j) \geq -x-1  \}}  e^{-\frac{1}{4}(x+V(w_{k-1}))}   \right) + 
\\
&\sum_{p\geq 0 } h(ce^{p+1}) \sum_{k=1}^n \E_\Q\left( R^-(x+1+V(w_k)) 1_{\{ \min_{j\in[0,k]}V(w_j) \geq -x-1,\, 4\ln \tilde{\Delta_k}-1  \leq  x+V(w_{k-1})  +p  \leq  4\ln \tilde{\Delta_k}   \}}    \right).
\end{align*}
The first term is bounded uniformly in $x\in \r^+$ and $n\in \N$ thanks to Lemma \ref{techniqAID}. For the second term let us introduce $(\tilde{\Delta},\zeta )$ be generic random variable distributed as $(\tilde{\Delta}_1,V(w_1)_+)$ under $\Q$. We can re-write this term as
\begin{align*}
\sum_{p\geq 0} h(ce^{p+1}) \sum_{k=1}^\infty & \E\left( R^-(4\ln \tilde{\Delta} -p+1 + \zeta) 1_{\{ \min_{j\in[0,k]}S_j \geq -x-1,\, 4\ln \tilde{\Delta}-1  \leq  x+S_{k-1}  +p  \leq  4\ln \tilde{\Delta}   \}}    \right)
\\
&\leq   \E\left( R^-(4\ln \tilde{\Delta} -p+1 + \zeta)  H(x+1, 4\ln \tilde{\Delta}-p)   \right)
\end{align*}
with $\forall x,r>0,\, 
H(x,r):= \sum_{k=1}^{+\infty} \P\left(  \min_{j\in[0,k-1]}S_j \geq -x,\, r-1  \leq  x+S_{k-1}  \leq  r    \right)$. Moreover by using the same arguments as in the proof of Lemma \ref{Lemconsequec}, it is plain to check that $ H(x,r)  \leq c (1+r) \1_{r \geq 0}$. Finally we get that
\begin{align*}
\sum_{k=1}^n I_k &\leq c+ \sum_{p\geq 0} h(e^{p+1}) \E_{\Q} \left(  R^-(4\ln \tilde{\Delta} + \zeta) 1_{\{ p -1 \leq  4\ln \tilde{\Delta}   \}}  ( 4\ln\tilde{\Delta}-p+1  )   \right)
\\
&\leq c+  \E_{\Q}\left(  R^-(4\ln \tilde{\Delta} + \zeta) \sum_{p=0}^{4\ln \tilde{\Delta_k}} (\ee-1 +\ln^2(e^{p+1})) ( 4\ln\tilde{\Delta}-p+1  )   \right) 
\\
&\leq c+ \E_{\Q}\left(  R^-(4\ln \tilde{\Delta} + \zeta) (\ln \tilde{\Delta})^4   \right) \leq  c \E((\tilde{X}+X)\max (0,\log \tilde{X}+X)^5)<\infty,
\end{align*}
where in the last line we used that $R^-(x)\leq c(1+x)$,  $\tilde{\Delta}+ \zeta$ is stochastically dominated by $2(\tilde{X}+X)$ and hypothesis \eqref{ourAssum}.\hfill$\Box$

  \begin{Lemma}
  	\label{toreste}
  Assume (\ref{chap3criticalcondition1}), (\ref{chap3criticalcondition2}) and  (\ref{ourAssum}). Let $\mathtt{u}\in \mathbb{T}$ be the vertex such that $V(u)=\M$ (if several such a vertex $\u$ exist one chooses one at random among the youngest one). For any $\epsilon >0$, 
  	\begin{align*}
  	\lim_{p\to\infty} \sup_{x\in \r^+} \P\left(\mathfrak{D}^\M  - \mathfrak{D}^{\u,\geq p}   \geq \epsilon \Big| \M\leq -x\right)=0.
  	\end{align*}	
  \end{Lemma}
  \noindent{\it Proof of Lemma \ref{toreste}.} The proof is quite similar to this one of Lemma \ref{cuty2}. According to Lemma \ref{tight}, it suffices to prove that 
  \begin{align*}
\lim_{p\to\infty}\sup_{x\in \r^+}e^x\P\left( \sum_{k=1}^{|\u|-p} \sum_{v \in \B(\u_k)} e^{\M-V(v) }   D_\infty^{(v)} \geq \epsilon ,\, \M\in I(x)\right) =0.
  \end{align*}
Recalling that $\1_{\{ \M \in I(x) \}} \leq   \sum_{k\geq 0} \sum_{|z|=k}  \1_{\{ V(z)=\M<\M_{|z|-1},\,    V(z)\in I(x) \}}$ we have
 \begin{align*}
 &e^x\P\left( \sum_{k=1}^{|\u|-p} \sum_{v \in \B(\u_k)} e^{\M-V(v) }   D_\infty^{(v)} \geq \epsilon ,\, \M\in I(x)\right) 
 \\
 &\leq e^x \E\left(  \sum_{n\geq p+1}  \sum_{|z|=n} \1_{\{ V(z)=\M<\M_{|z|-1},\,    V(z)\in I(x) \}}   1_{\{ \sum_{k=1}^{n-p} \sum_{v \in \B(z_k)} e^{\M-V(v) }   D_\infty^{(v)} \geq \epsilon \}}  \right)
 \\
 &\overset{\text{Prop \ref{lyons}}}{\leq } e^x \sum_{n\geq p+1} \E\left(    e^{V(w_n)} \1_{\{ V(w_n)<\min_{j\in [0,n-1]}V(w_j),\,    V(w_n)\in I(x) \}}   1_{\{ \sum_{k=1}^{n-p} \sum_{v \in \B(w_k)} e^{\M-V(v) }   \bar{D}_\infty^{(v,x+1)} \geq \epsilon \}}  \right).
 \end{align*}
 where we recall that $\bar{D}_\infty^{(v,x+1)}:= \lim_{n\to\infty} \sum_{|z|=n,\, z\geq v} [V(z)-V(v)]e^{-V(z)-V(v)} \1_{\{ \min_{i\in [|v|,n]}V(z_i)\geq -x-1-V(v)  \}}$. By conditioning with respect to the sigma-field $\sigma(\B_n):= \sigma (  V(w_j), (V(u),\, u\in \B(w_j)))_{j\in [1,n]}$, we get

 \begin{align*}
& e^x\P\left( \sum_{k=1}^{|\u|-p} \sum_{v \in \B(\u_k)} e^{\M-V(v) }   D_\infty^{(v)} \geq \epsilon ,\, \M\in I(x)\right) 
 \\
  &\leq  e^x \sum_{n\geq p+1} \E\left(    e^{V(w_n)} \1_{\{ V(w_n)<\min_{j\in [0,n-1]}V(w_j),\,    V(w_n)\in I(x) \}}   \P\left( \sum_{k=1}^{n-p} \sum_{v \in \B(w_k)} e^{\M-V(v) }   \bar{D}_\infty^{(v,x+1)} \geq \epsilon \Big| \sigma(\B_n) \right)  \right)
 \\
 &\leq \sum_{n\geq p+1} \E\left(     \1_{\{ V(w_n)<\min_{j\in [0,n-1]}V(w_j),\,    V(w_n)\in I(x) \}}  \min(1, \frac{1}{\epsilon} \sum_{k=1}^{n-p}    \sum_{v \in \B(w_k)} (x+1+V(v)) e^{-x-V(v) } )    \right),
 \end{align*}
 where in the last line we used the Markov inequality then Eq.(\ref{piousit}). Let\\ $\tilde{\Delta}_k:=   \sum_{v \in \B(w_k)} e^{-[V(v)-V(w_{k-1})]} ( [V(v)-V(w_{k-1})]_+ +1)$, then notice that for any $k\in [1,n]$,  
 \begin{align*}
  \sum_{v \in \B(w_k)} (x+1+V(v)) e^{-x-V(v) }&\leq ([x+ V(w_{k-1})]_+ +1)e^{-x-V(w_{k-1})} \tilde{\Delta}_k
  \\
  &\leq3e^{-\frac{1}{2}(x+V(w_{k-1}))} \tilde{\Delta}_k.
 \end{align*}
On the other hand we have
 \begin{align*}
 \min(1, \frac{1}{\epsilon} \sum_{k=1}^{n-p}   3e^{-\frac{1}{2}(x+V(w_{k-1}))} \tilde{\Delta}_k ) \leq  \frac{3}{\epsilon}e^{-\frac{1}{4}(x+V(w_{k-1}))} + 1_{\{ x+V(w_{k-1}) \leq 4\ln \tilde{\Delta}_k   \}}.
 \end{align*}
Thus it remains to prove the following two limits
 \begin{align}
\label{visit1} \lim_{p\to\infty}\sup_{x\in \r^+ }  \sum_{n\geq p+1} \sum_{k=1}^{n-p}\E\left(      \1_{\{ S_n<\min_{j\in [0,n-1]}S_j,\,    S_n\in I(x) \}}       e^{-\frac{1}{4}(x+S_k) }     \right)=0,
\\
\label{visit2}  \lim_{p\to\infty}\sup_{x\in \r^+} \sum_{n\geq p+1} \sum_{k=1}^{n-p} \P\left(  V(w_n)<\min_{j\in [0,n-1]}V(w_j),\,    V(w_n)\in I(x) ,\, x+V(w_{k-1}) \leq 4\ln \tilde{\Delta}_k       \right)=0.
 \end{align}
 By operating a time reversal on the random walk $(S_j)_{j\leq n}$ we have
 \begin{align*}
 \sum_{n\geq p+1} \sum_{k=1}^{n-p} \E\left(      \1_{\{ S_n<\min_{j\in [0,n-1]}S_j,\,    S_n\in I(x) \}}       e^{-\frac{1}{4}(x+S_{k-1}) }    \right)    \leq  \sum_{n\geq p+1} \sum_{k=1}^{n-p} \E\left(    \1_{\{ \max_{j\in [1,n]}S_j<0,\,   S_n\in I(x) \}}      e^{\frac{1}{2}S_{n-k+1} }    \right)\qquad\qquad 
\\
=  \sum_{k=p+1}^{+\infty} \E\left(    \1_{\{ \max_{j\in [1,k]}S_j<0\}} e^{\frac{1}{2}S_k} \sum_{n\geq 0} \P_{S_k}\left(  \max_{j\in [1,n]}S_j <0   ,\,    S_{n}\in I(x) \right)_{z=S_k}        \right)
\\
= \sum_{k=p+1}^{+\infty} \E\left(    \1_{\{ \max_{j\in [1,k]}S_j<0\}} e^{\frac{1}{2}S_k} \sum_{n\geq 0} [ \tilde{R}(x+1,-S_k) -\tilde{R}(x,-S_k) ]     \right),\qquad\,\,
\end{align*}
where in the last line we  inverted the sums and used the Markov property. Finally by using Lemma \ref{Lemconsequec} and Lemma \ref{techniqAID} we deduce that for any $x>0$,
\begin{align*}
\sum_{n\geq p+1} \sum_{k=1}^{n-p}\E\left(      \1_{\{ S_n<\min_{j\in [0,n-1]}S_j,\,    S_n\in I(x) \}}       e^{-\frac{1}{4}(x+S_j) }     \right) &\leq c \sum_{k=p+1}^{+\infty} \E\left(    \1_{\{ \max_{j\in [1,k]}S_j<0\}} e^{\frac{1}{2}S_k}  (1-S_k)    \right) \underset{p\to\infty}{\to} 0,
 \end{align*}
which proves (\ref{visit1}). It remains to prove (\ref{visit2}). Recall that $(V(w_j)-V(w_{j-1}),\tilde{\Delta}_j)_{j\in [|1,k|]}$ are i.i.d., then by operating a time reversal we get for any $x\in \r^+$, 
 \begin{align*}
 &\sum_{n\geq p+1} \sum_{k=1}^{n-p} \P\left(  V(w_n)<\min_{j\in [0,n-1]}V(w_j),\,    V(w_n)\in I(x) ,\, x+V(w_{k-1}) \leq 4\ln \tilde{\Delta}_k       \right)
 \\
& =   \sum_{n\geq p+1} \sum_{k=1}^{n-p} \P\left(  \max_{j\in [1,n]}V(w_j)<0,\,    V(w_n)\in I(x) ,\, V(w_n) +x - V(w_{n-k+1}) \leq 4\ln \tilde{\Delta}_{n-k+1}    \right)
 \\
 &\leq  \sum_{n\geq p+1} \sum_{k= p+1}^{n} \P\left(  \max_{j\in [1,n]}V(w_j)<0,\,    V(w_n)\in I(x) ,\, - V(w_{k}) \leq 4\ln \tilde{\Delta}_{k}  +1   \right).
 \end{align*}
 By using the Markov property at time $k$, it follows that 
 \begin{align*}
 &\sum_{n\geq p+1} \sum_{k=1}^{n-p} \P\left(  V(w_n)<\min_{j\in [0,n-1]}V(w_j),\,    V(w_n)\in I(x) ,\, x+V(w_{k-1}) \leq 4\ln \tilde{\Delta}_k       \right)
 \\
  &=  \sum_{k= p+1}^{+\infty} \E\left( 1_{\{ \max_{j\in [1,k]}V(w_j)<0,\, - V(w_{k}) \leq 4\ln \tilde{\Delta}_{k}  +1 \}} \sum_{n\geq 0} \E_{V(w_k)}\left(  \max_{j\in [1,n]}S_j<0 ,\,S_n\in I(x)  \right)  \right)
 \\
 &= \sum_{k= p+1}^{+\infty} \E\left( 1_{\{ \max_{j\in [1,k]}V(w_j)<0,\, - V(w_{k}) \leq 4\ln \tilde{\Delta}_{k}  +1 \}}  [ \tilde{R}(x+1,-V(w_k)) - \tilde{R}(x,-V(w_k)) ]\right)
 \\
 &\leq c\sum_{k= p+1}^{+\infty} \E\left( 1_{\{ \max_{j\in [1,k]}V(w_j)<0,\, - V(w_{k}) \leq 4\ln \tilde{\Delta}_{k}  +1 \}} (1+ 4\ln \tilde{\Delta}_{k}  +1 )\right),
 \end{align*}
where we used Lemma \ref{Lemconsequec} in the last inequality. By introducing $(\tilde{\Delta},\zeta)$ a random variable independent of everything and distributed as $(\tilde{\Delta}_1,V(w_1))$, we get that for any $x\in \r, p>0$, 
\begin{align*}
&\sum_{n\geq p+1} \sum_{k=1}^{n-p} \P\left(  V(w_n)<\min_{j\in [0,n-1]}V(w_j),\,    V(w_n)\in I(x) ,\, x+V(w_{k-1}) \leq 4\ln \tilde{\Delta}_k       \right)
\\
&\leq c \E\left(  (1+ 4\ln \tilde{\Delta}  +1 )\sum_{k= p+1}^{+\infty}\P\left( \max_{j\in [1,k]}S_j<0,\, - S_{k-1} \leq 4\ln \tilde{\Delta}  +1 +\zeta \right) \right)
\\
& \leq c' \E \left( (\tilde{X}+X)\max (0,\ln \tilde{X}+X)^4 \right) <+\infty,
\end{align*}
where we used (\ref{ourAssum}) in the last inequality. The sum in the second line does not depend in $x$ any more and is finite, thus when $p$ goes to $\infty$ the sum converges toward zero which concludes the proof of Lemma \ref{toreste}.\hfill$\Box$

 \bibliographystyle{plain}
 \bibliography{bibli}
 
\end{document}